\documentstyle[11pt]{amsart}

\newcommand{\nc}{\newcommand}

\nc{\ad}{\operatorname{ad}}
\nc{\Tot}{\operatorname{Tot}}
\nc{\Aut}{\operatorname{Aut}}
\nc{\boxtimes}{\fbox{$\times$}}
\nc{\CCF}{\cal{CF}}
\nc{\CDer}{\cal{D}er}
\nc{\CDF}{\cal{DF}}
\nc{\CHC}{\check{\cal C}}
\nc{\CHom}{\cal{H}om}
\nc{\Cone}{\operatorname{Cone}}
\nc{\dec}{\operatorname{dec}}
\nc{\Der}{\cal{D}er}
\nc{\Diff}{\cal{D}\mbox{\it iff}}
\nc{\dirlim}{\underset{\rightarrow}{\operatorname{lim}}}
\nc{\dpar}{\partial}
\nc{\GL}{\operatorname{GL}}
\nc{\CGr}{\cal{G}r}
\nc{\gr}{\operatorname{gr}}
\nc{\pr}{\operatorname{pr}}
\nc{\semid}{|\!\!\!\times}
\nc{\Hom}{\operatorname{Hom}}
\nc{\Id}{\operatorname{Id}}
\nc{\id}{\operatorname{id}}
\nc{\Ima}{\operatorname{Im}}
\nc{\invtimes}{\underset{\leftarrow}{\otimes}}
\nc{\invlim}{\underset{\leftarrow}{\operatorname{lim}}}
\nc{\Lie}{\operatorname{Lie}}
\nc{\res}{\operatorname{res}}
\nc{\Spec}{\operatorname{Spec}}
\nc{\Spf}{\operatorname{Spf}}
\nc{\sym}{\operatorname{sym}}
\nc{\Td}{\operatorname{Td}}
\nc{\Th}{\mbox{\bf T}}
\nc{\Tor}{\operatorname{Tor}}

\nc{\Artin}{\cal{A}rtin}
\nc{\Dgcoalg}{\cal{D}gc}
\nc{\Dgcu}{\cal{D}gc^u}
\nc{\Dglie}{\cal{D}glie}
\nc{\Ens}{\cal{E}ns}
\nc{\Fsch}{\cal{F}sch}
\nc{\Gr}{\cal{G}r}
\nc{\Groupoids}{\cal{G}roupoids}
\nc{\Holie}{\cal{H}olie}
\nc{\Mod}{\cal{M}od}
\nc{\Mor}{\cal{M}or}
\nc{\Sch}{\cal{S}ch}

\nc{\CA}{\cal A}
\nc{\CB}{\cal B}
\nc{\CC}{\cal C}
\nc{\CD}{\cal D}
\nc{\CF}{\cal F}
\nc{\CG}{\cal G}
\nc{\CH}{\cal H}
\nc{\CI}{\cal I}
\nc{\CM}{\cal M}
\nc{\CN}{\cal N}
\nc{\CO}{\cal O}
\nc{\CP}{\cal P}
\nc{\CT}{\cal T}
\nc{\CU}{\cal U}
\nc{\CV}{\cal V}
\nc{\CZ}{\cal Z}

\nc{\fa}{\frak{a}}
\nc{\fA}{\frak{A}}
\nc{\fg}{\frak{g}}
\nc{\fh}{\frak{h}}
\nc{\fI}{\frak{I}}
\nc{\fK}{\frak{K}}
\nc{\fm}{\frak{m}}
\nc{\fP}{\frak{P}}
\nc{\fS}{\frak{S}}
\nc{\ft}{\frak{t}}
\nc{\fX}{\frak{X}}
\nc{\fY}{\frak{Y}}

\nc{\bF}{\bar{F}}
\nc{\bCP}{\bar{\cal{P}}}
\nc{\bm}{\mbox{\bf{m}}}
\nc{\bT}{\mbox{\bf{T}}}
\nc{\hP}{\hat{P}}

\nc{\nen}{\newenvironment}
\nc{\ol}{\overline}
\nc{\ul}{\underline}
\nc{\ra}{\rightarrow}
\nc{\lla}{\longleftarrow}
\nc{\lra}{\longrightarrow}
\nc{\Lra}{\Longrightarrow}
\nc{\Lla}{\Longleftarrow}
\nc{\Llra}{\Longleftrightarrow}
\nc{\hra}{\hookrightarrow}
\nc{\iso}{\overset{\sim}{\lra}}

\nc{\Thm}[1]{Theorem~\ref{#1}}
\nc{\Prop}[1]{Proposition~\ref{#1}}
\nc{\Lem}[1]{Lemma~\ref{#1}}
\nc{\Cor}[1]{Corollary~\ref{#1}}
\nc{\Conj}[1]{Conjecture~\ref{#1}}
\nc{\Claim}[1]{Claim~\ref{#1}}
\nc{\Defn}[1]{Definition~\ref{#1}}
\nc{\Exa}[1]{Example~\ref{#1}}
\nc{\Rem}[1]{Remark~\ref{#1}}
\nc{\Note}[1]{Note~\ref{#1}}

\nen{thm}[1]{\label{#1}{\bf Theorem.\ } \em}{}
\nen{prop}[1]{\label{#1}{\bf Proposition.\ } \em}{}
\nen{lem}[1]{\label{#1}{\bf Lemma.\ } \em}{}
\nen{cor}[1]{\label{#1}{\bf Corollary.\ } \em}{}
\nen{conj}[1]{\label{#1}{\bf Conjecture.\ } \em}{}
\nen{claim}[1]{\label{#1}{\bf Claim.\ } \em}{}

\nen{defn}[1]{\label{#1}{\bf Definition.\ } }{}
\nen{exa}[1]{\label{#1}{\bf Example.\ } }{}

\nen{rem}[1]{\label{#1}{\em Remark.\ } }{}
\nen{note}[1]{\label{#1}{\em Note.\ } }{}
\nen{exer}[1]{\label{#1}{\em Exercise.\ } }{}

\begin{document}



\title[]{Deformation theory and Lie algebra homology}
\author{Vladimir Hinich}
\address{V.H.: Dept. of Mathematics and Computer Science, University of Haifa,
Mount Carmel, Haifa 31905 Israel}
\author{Vadim Schechtman}
\address{V.Sch.: Dept. of Mathematics, SUNY at Stony Brook, Stony Brook
NY 11794 USA}

\date{May 1994}

\thanks{The work was partially supported by the Paul and Gabriella Rosenbaum
Foundation. The second author was supported in part by the NSF grant
DMS-9202280.}

\maketitle


\section{Introduction}

\subsection{}
\label{pose} Let $X$ be a smooth proper scheme $X$ over a field $k$
of characteristic $0$, $G$ an algebraic group over $k$ and $p: P\lra X$
a $G$-torsor over $X$. Consider the following deformation
problems.

{\bf Problem 1.} Flat deformations of $X$.

{\bf Problem 2.} Flat deformations of the pair $(X,P)$.

{\bf Problem 3.} Deformations of $P$ ($X$ being fixed).

According to Grothendieck, one  can assign to a deformation problem
of the above type a sheaf of Lie algebras
over $X$ which can be defined as "a sheaf of infinitesimal automorphisms"
of the corresponding deformation functor (cf. Section ~\ref{univers}).

Let us describe the sheaves corresponding to our problems. For Problem 1 it is
a tangent sheaf $\CA_1=\CT_X$. For Problem 2 it is a sheaf $\CA_2=\CA_P$
defined as follows. For a Zariski open $U\subset X$, $\Gamma(U,\CA_P)$
is the space
of $G$-invariant vector fields on $p^{-1}(U)$. The map $p$ induces
(epimorhic) map $\epsilon:\CA_P\lra\CT_X$. For Problem 3, $\CA_3=\fg_P:=
\ker(\epsilon)$. The last sheaf may be also defined as a sheaf of sections
of a vector bundle associated with $P$ and the adjoint representation
of $G$ on its Lie algebra. Note that the sheaves $\CA_i$ are locally free
$\CO_X$-modules of finite type (but the bracket is not $\CO_X$-linear for
$i=1,2$).

Suppose that $H^0(X,\CA_i)=0$. Then it is known that there exists
a universal formal deformation space $\fS_i=\Spf(R_i)$ for Problem $i$
(see Section ~\ref{univers}). Here $R_i$ is a complete local $k$-algebra
with residue field $k$. Let $\bm_{R_i}$ denote its maximal ideal.
We have the Kodaira-Spencer isomorphism
\begin{equation}
\label{ksg}
\kappa^1: T_{\fS_i,s}=(\bm_{R_i}/\bm_{R_i}^2)^*\iso H^1(X,\CA_i)
\end{equation}
describing the tangent space of $\fS_i$ at the closed point $s$.

The  main goal of the present paper is to describe, in case when $\fS_i$ is
smooth, the whole ring $R_i$ in terms of the sheaf $\CA_i$.

\subsection{}
\label{sullivan} To formulate the answer, we need a certain cohomological
construction. Let $\fg$ be a sheaf of $k$-Lie algebras over $X$ which is
also a quasicoherent
$\CO_X$-module. Consider an affine open covering $\CU$
of $X$ and the corresponding complex of \v{C}ech cochains $\CHC(\CU,\fg)$.
Using a generalization of Thom-Sullivan construction
\footnote{introduced in ~\cite{hlha}} used in Rational
homotopy theory one can construct (see Section ~\ref{direct})
a certain differential graded Lie algebra
$R\Gamma^{Lie}(X,\fg)$ canonically quasi-isomorphic to $\CHC(\CU,\fg)$.
This dg Lie algebra does not depend, up to (essentially) unique
quasi-isomorphism, on a covering, hence we omit it from the notation.

Now one can apply to $R\Gamma^{Lie}(X,\fg)$ the Quillen functor $C$
(which is a generalization to dg Lie algebras of the Chevalley complex
computing the homology of a Lie algebra with trivial coefficients) and
get the complex $C(R\Gamma^{Lie}(X,\fg))$. This complex carries the canonical
increasing filtration $\{ F_nC(R\Gamma^{Lie}(X,\fg))\}$ with
graded quotients isomorphic to symmetric powers
$S^n(R\Gamma^{Lie}(X,\fg)[1])$ (for details, see \ref{quillen}).
Let us define homology spaces
$$
H^{Lie}_i(R\Gamma^{Lie}(X,\fg)):=H^{-i}(C(R\Gamma^{Lie}(X,\fg)));\
F_nH^{Lie}_i(R\Gamma^{Lie}(X,\fg)):=H^{-i}(F_nC(R\Gamma^{Lie}(X,\fg))).
$$
These spaces depend only on $X$ and on the sheaf $\fg$.
We have evident maps
$$
\ldots\lra F_{n-1}H^{Lie}_i(R\Gamma^{Lie}(X,\fg))\lra
F_nH^{Lie}_i(R\Gamma^{Lie}(X,\fg))\lra\ldots\lra
H^{Lie}_i(R\Gamma^{Lie}(X,\fg))
$$
and
$$
F_nH^{Lie}_i(R\Gamma^{Lie}(X,\fg))\lra
H^{n-i}(\Lambda^iR\Gamma^{Lie}(X,\fg))
$$
($\Lambda^i$ denotes the exterior power).
\subsection{} Return to the assumptions ~\ref{pose}. Pick $i=1,2$ or $3$.
Recall that we suppose that $H^0(X,\CA_i)=0$. Suppose also that
$\fS_i$ is {\em smooth} i.e. isomorphic to a formal power series ring
over $k$.

Let $R^*_i$ denote the space of continuous $k$-linear maps
$R_i\lra k$ ($k$ is equipped with the discrete topology). The main result
of this paper is (see Thm.~\ref{complet}).

\subsubsection{}
\label{thm-intro}
\begin{thm}{} One has compatible canonical isomorphisms
\begin{equation}
\label{isomor}
\kappa: R^*_i\iso H^{Lie}_0(R\Gamma^{Lie}(X,\CA_i));\
\end{equation}
\begin{equation}
\kappa^{\leq n}: (R_i/\bm^{n+1}_{R_i})^*\iso
F_nH^{Lie}_0(R\Gamma^{Lie}(X,\CA_i))
\end{equation}
After passing to the graded quotients, isomorphisms $\kappa^{\leq n}$ induce
isomorphisms
$$
S^n(T_{\fS_i,s})\cong (\bm^{n}/\bm^{n+1})^*\iso
H^{n}(\Lambda^nR\Gamma^{Lie}(X,\CA_i))\cong S^nH^1(X,\CA_i)
$$
which coincide with $(-1)^n\kappa^1$.
\end{thm}

\subsection{} Let us describe the contents of the paper in a more detail.

Our construction of the isomorphism $\kappa$ is based on the construction of
{\bf higher Kodaira-Spencer morphisms} --- they are relative versions
of ~(\ref{isomor}). Similarly to a usual KS map they are not necessarily
isomorphisms and are defined for arbitrary --- not necessarily universal ---
deformation\footnote{this idea is used in ~\cite{bs}}.

In the same way as the usual KS map comes from a certain boundary homomorphism
of a short exact sequence of complexes, our higher KS map comes from
a map, which we call {\bf connecting morhism}. It arises from
an extension of {\em dg Lie algebras}
$$
0\lra\fg\lra\fa\lra\ft\lra 0
$$
and maps the enveloping algebra of the cone of the map $\fg\ra\fa$
(which is quasi-isomorphic to $U(\ft)$) to the standard complex
$C(\fg)$.
The construction is given in Section ~\ref{envstand}, cf. ~\ref{constr-thm}.

For applications we need a generalization of this construction to
{\em dg Lie algebroids}. The simplest example of a sheaf of Lie algebroids is
the tangent sheaf $\CT_S$ over a variety $S$. We borrow
their definition from ~\cite{bb} and ~\cite{bfm}. Following
{\em loc.cit.}, one defines an
analogue of enveloping algebra for Lie algebroids --- we call them
{\em twisted enveloping algebras}. For example, the twisted enveloping algebra
of $\CT_S$ will be the sheaf $\Diff_S$ of differential operators
provided $S$ is smooth.
These sheaves are sources of higher KS maps. Section ~\ref{twisted}
generalizes the principal construction of Section ~\ref{envstand} to
(sheaves of) dg Lie algebroids.

The construction of connecting morphisms
is our first point. (They may be of an independent interest.)

Section ~\ref{formal} is technical; we discuss there some basic
properties of differential operators on formal schemes.

\subsection{} In Sections~\ref{direct} and~\ref{thoms} we develop
certain formalism of {\bf Homotopy Lie algebras}. We define there for
(locally noetherian)
schemes $X$ categories $\Holie^{qc}(X)$ which contain sheaves
of dg Lie algebras with bounded cohomology and flat quasicoherent components,
and have good functorial properties. For example, the functor $R\Gamma^{Lie}$
mentioned in ~\ref{sullivan} is a particular case of the {\bf direct image}
functor --- its construction is the main result of these Sections.

Although it is a technical tool, we regard the contents of Sect.~\ref{direct}
and ~\ref{thoms}  as a {\bf second main outcome} of the present paper.

We believe that unlike usual types of algebras,
{\em Homotopy Lie algebras} are correctly defined only as objects of
appropriate {\em Homotopy categories}. The main property one needs
from them is good functorial behaviour. This is what Thom-Sullivan
functor does. The success of the construction is due to the fact that
Thom-Sullivan functor has excellent exactness, base change, etc.
properties. They are proved in ~\ref{thoms}.

Another approach to the definition and functoriality of Homotopy Lie algebras
has been developed in ~\cite{hla}. In {\em loc.cit.} a HLA was a
"dg Lie algebra up to higher homotopies". There are strong indications that
both definitions give rise to equivalent Homotopy categories.

It seems that one can use also cosimplicial Lie algebras
for the construction of the same Homotopy category.
 From this point of view, different versions of HLA are nothing but
different "models" of these Homotopy categories. We believe that
every "model" might be useful in applications.

Of course these remarks apply also to other types of algebras.

\subsection{} In Section ~\ref{kodaira}  we take all
the results together and start the machine of Section ~\ref{twisted} which
cooks up higher KS maps for us. In Section ~\ref{univers} we discuss
universal formal deformations.

\subsection{} One can generalize the above approach and get a description
of sheaves of differential operators on moduli spaces acting on natural
vector bundles (for example, determinant bundles) in terms of Lie algebra
homology (with non-trivial coefficients).
For example, for operators of order $\leq 1$ on determinant
bundles one gets a result equivalent to ~\cite{bs}.

Also there are strong indications that combining our description with
Serre duality it is possible to get a "global counterpart" of the results
of ~\cite{bd}.

We will return to these subjects in future publications.

\subsection{} The general main idea that

{\bf the completion of a local ring of a moduli
space at a given point $X$ is isomorphic to the dual of the $0$-th
homology group
of the "Lie algebra of infinitesimal automorphisms of $X$"}

was spelled out very clearly a few years ago by Drinfeld,
Deligne, Feigin  (cf. ~\cite{d}, ~\cite{del}, \cite{f}). We knew this
idea from Drinfeld.

For deformations of complex structures a result analogous
to Theorem ~\ref{thm-intro} was proven in ~\cite{gm} (by different argument).

The present paper develops further the results of ~\cite{hdt}.
Note also some recent related work: ~\cite{ev}, \cite{r}.

This paper owes much to various ideas of A.Beilinson and V.Drinfeld.
We express to them our deep gratitude. We thank H.~Esnault, V.~Ginzburg
and E.~Viehweg for useful discussions.

We are especially grateful to Professor Han Sah who made possible several
visits of V.H. to Stony Brook.

\section{Enveloping algebras and standard complexes}
\label{envstand}

\subsection{Preliminaries}

\subsubsection{} Let us fix some notations and sign conventions
(cf. ~\cite{de}, 1.1).

Throughout this Section we fix a commutative ground ring $k$ of containing
$\Bbb Q$.

$\Mod (k),\ \Gr (k)$ will denote a category of $k$-modules and that of
$\Bbb Z$-graded $k$-modules respectively. $\CC (k)$
denotes the category of complexes of $k$-modules (all differentials have
degree $+1$). We have an obvious forgetful functor $\CC (k)\lra \Gr (k)$.
If $X\in \CC (k)$ (or $\Gr (k)$), $x\in X^p$, we refer to $p$ as to degree of
$x$, and denote it $|x|$.

We identify $k$-modules with complexes concentrated in degree $0$.
If $n\in \Bbb Z$, $X\in \Gr (k)$,  $X[n]$ will denote the shifted module
$X[n]^p=X^{p+n}$. If $X\in \CC(k)$, we define $X[n]\in \CC(k)$ which is as
above as a graded module, the differential being $d_{X[n]}=(-1)^nd_X$.

{\em A map of degree $n$}, $f: X\lra Y$ in $\CC(k)$ is by definition a
morphism $f: X\lra Y[n]$ in $\Gr (k)$. For such an $f$ we set $df=d_Y\circ f+
(-1)^nf\circ d_X: X\lra Y[n+1]$.

The category $\CC (k)$ has a tensor structure --- for $X,Y\in \CC (k)$,
$X\otimes Y$ is the usual tensor product of complexes over $k$. We have
natural associativity and commutativity isomorphisms
\begin{equation}
\label{assoc}
a_{X,Y,Z}:(X\otimes Y)\otimes Z\iso X\otimes (Y\otimes Z),\
(x\otimes y)\otimes z\mapsto x\otimes (y\otimes z)
\end{equation}
and
\begin{equation}
\label{r}
R_{X,Y}: X\otimes Y\iso Y\otimes X
\end{equation}
defined by the formula $R_{X,Y}(x\otimes y)=(-1)^{|x||y|}y\otimes x$.

This formula, as well as other formulas in the "dg world" are obtained by
implementing {\em  the Quillen sign rule}:
"whenever something of degree
$p$ is moved past something of degree $q$ the sign $(-1)^{pq}$ accrues",
{}~\cite{q}, p. 209.

Isomorphisms ~(\ref{assoc}) and ~(\ref{r}) endow $\CC(k)$ with a structure
of a
{\em symmetric monoidal category} in the sense of MacLane, ~\cite{mac}.

We have canonical {\em shifting} isomorphisms
\begin{equation}
\label{shift}
X[n]\otimes Y[m]\iso (X\otimes Y)[n+m]
\end{equation}
sending $x\otimes y$, $x\in X^i,\ y\in Y^j$ to $(-1)^{im}x\otimes y$.
(They may be obtained by identifying (following ~\cite{de}, 1.1) $X[n]$ with
$k[n]\otimes X$ , and applying $R_{X,k[m]}$).

{\em Algebras.} The structure of a symmetric monoidal category allows one
to define usual types of algebras in $\CC (k)$.
We shall refer to them by adding "dg" to their original name.
More specifically, we shall use the following algebras (cf. ~\cite{q}).

--- A {\em dg Lie algebra} is a complex $X$ together with a bracket
$[,]: X\otimes X\lra X$ which is skew symmetric, i.e.
$[,]\circ R_{X,X}=-[,]$, and satisfies the Jacobi identity
$$
[x,[y,z]]+(-1)^{|x|(|y|+|z|)}[y,[z,x]]+(-1)^{|z|(|x|+|y|)}[z,[x,y]]=0
$$
We shall denote $\Dglie=\Dglie_k$ the category of dg Lie algebras
over $k$.

--- A {\em dg coalgebra} is a complex $X$ together with an coassociative
comultiplication $\Delta: X\lra X\otimes X$ and a counit $\epsilon: X\lra k$.
$X$ is called {\em cocommutative} if $\Delta= R_{X,X}\circ \Delta$.
We shall denote $\Dgcoalg=\Dgcoalg_k$ the category of cocommutative
dg coalgebras over $k$.

--- A {\em dg algebra} is a complex $X$ together with an associative
multiplication $\mu: X\otimes X\lra X$ and the unit $1\in X^0$ (which
may be considered as a map $k\lra X$). $X$ is called {\em commutative}
is $\mu=\mu\circ R_{X,X}$.

--- A {\em dg Hopf algebra} is a complex $X$ together with a multiplication
$\mu$, comultiplication $\Delta$, a unit and a counit, making it
a dg algebra and a dg coalgebra, and such that these two structures are
compatible in the standard way. In particular, $\Delta$ is
a map of dg algebras, where the multiplication in $X\otimes X$ is defined
by the rule $(x\otimes y)(x'\otimes y')=(-1)^{|y||x'|}xx'\otimes yy'$.
$X$ is called (co)commutative if the underlying (co)algebra is.
An element $x\in X$ is called {\em primitive} if
$\Delta (x)=x\otimes 1+1\otimes x$.

\subsubsection{}
\label{connected}
An element $u$ of a dg coalgebra $C$ is called {\em group-like}
if the following conditions are fulfilled:

(i) $du=0$ (ii) $\Delta(u)=u\otimes u$ (iii) $\epsilon(u)=1\in k$

A group-like element $u\in C$ defines a decomposition $C=k\cdot u\oplus
C^+$ with $C^+=\ker(\epsilon)$.
Let $\pi_u:C\ra C^+$ be the projection onto the second summand.

For a positive integer $n$ define the map
$$
\Delta^{(n)}: C\lra C^{\otimes n}
$$
by induction on $n$: set $\Delta^{(1)}=\id_C,
\ \Delta^{(n)}=(\Delta\otimes \id_{C^{\otimes(n-2)}})\circ \Delta^{(n-1)}$.

The choice of a group-like element $u$ defines an increasing filtration
$F^u=\{F^u_n\}$ of $C$ by the formula

$$ F^u_n=\ker\left(C\overset{\Delta^{(n+1)}}{\lra}C^{\otimes n+1}
\overset{\pi_u^{\otimes n+1}}{\lra}(C^+)^{\otimes n+1}\right).$$

\begin{defn}{unital}
1. A group-like element $u\in C$ is called {\em a unit} if the corresponding
filtration $F^u$ is exhaustive: $C=\cup_{i=0}^{\infty} F_i$.

2. {\em A unital coalgebra} is a pair $(C,1_C)$ where $C$ is a coalgebra and
$1_C\in C$ is a unit it in.
\end{defn}

The category $\Dgcu(k)$ has as objects unital cocommutative dg $k$-coalgebras;
a morphism in $\Dgcu(k)$  is a coalgebra morphism preserving the units.

\subsubsection{} Let $X\in \CC (k)$. For an integer $n\geq 1$ denote
$T^n(X)$ its $n$-th tensor power $X\otimes\ldots \otimes X$. Set $T^0X=k$.
The direct sum $T(X)=\oplus_{n\geq 0} T^n(X)$ has a natural structure
of a dg algebra --- the {\em tensor algebra} of $X$.

The commutativity isomorphisms ~(\ref{r}) define the action of the symmetric
group on $n$ letters $\Sigma_n$ on $T^n(X)$. We denote $S^n(X)$, $\Lambda^n(X)$
the complexes of coinvariants (resp., coantiinvariants) of this action
and $\pi_{S,n}: T^n(X)\lra S^n(X)$,
$\pi_{\wedge,n}: T^n(X)\lra \Lambda^n(X)$
the canonical projections.

For $x_1,\ldots, x_n\in X$ we set
$$
x_1\cdot\ldots\cdot x_n=\pi_{S,n}(x_1\otimes\ldots \otimes x_n)\in S^n(X);
$$
$$
x_1\wedge\ldots\wedge x_n=\pi_{\wedge,n}(x_1\otimes\ldots \otimes x_n)\in
\Lambda^n(X)
$$

The projections $\pi_{S,n},\ \pi_{\wedge,n}$ have
canonical sections
\begin{equation}
\label{sym}
i_{S,n}: S^n(X)\lra T^n(X),\ i_{S,n}(x_1\cdot\ldots\cdot x_n)=
\frac{1}{n!}\sum_{\sigma\in \Sigma_n}\sigma(x_1\otimes\ldots\otimes x_n)
\end{equation}
\begin{equation}
\label{antisym}
i_{\wedge,n}: \Lambda^n(X)\lra T^n(X),\ i_{\wedge,n}(x_1\cdot\ldots\cdot x_n)=
\frac{1}{n!}\sum_{\sigma\in \Sigma_n}(-1)^{|\sigma|}
\sigma(x_1\otimes\ldots\otimes x_n)
\end{equation}

The isomorphisms ~(\ref{shift}) induce  canonical isomorphisms
$$
a_n: T^n(X[1])\cong T^n(X)[n]
$$
such that  for $\sigma\in \Sigma_n$, $a_n\circ \sigma=(-1)^{|\sigma|}\sigma
\circ a_n$. After passing to coinvariants, we get canonical shifting
({\em "d\'{e}calage"}) isomorphisms
\begin{equation}
\label{dec}
\dec_n:S^n(X[1])\cong \Lambda^n(X)[n]
\end{equation}

Explicit formula: for $x_i\in X^{p_i},\ i=1,\ldots , n$,
\begin{equation}
\label{decalage}
\dec_n(x_1\cdot\ldots\cdot x_n)=(-1)^{\sum_{i=1}^n (n-i)p_i}
(x_1\wedge\ldots\wedge x_n)
\end{equation}

\subsubsection{}
The direct sum $S(X)=\sum_{n\geq 0} S^n(X)$ is naturally a commutative
dg algebra. Let us define the counit $S(X)\lra k=S^0(X)$ as a canonical
projection, and a map $\Delta: S(X)\lra S(X)\otimes S(X)$ by two conditions
which characterize it uniquely:

(i) $\Delta(x)=x\otimes 1 + 1\otimes x$ for $x\in X=S^1(X)$

(ii) $\Delta$ is a map of dg algebras.

This makes $S(X)$ a commutative and cocommutative dg Hopf algebra, cf.
{}~\cite{q},
App B, 3.3.

The argument of {\em loc.cit.} shows that $S(X)$ is a unital dg coalgebra
with
$1_{S(X)}=1\in k=S^0(X)$. The corresponding filtration is
$F_iS(X)=\oplus_{p=0}^iS^p(X)$.

\subsubsection{Universal property of $S(X)$}
\label{coalg-s}

Let $C\in\Dgcu$  and suppose we are given a map of unital dg coalgebras
$f: C\lra S(X)$. Set $f_n=p_n\circ f: C\lra S^n(X)$ where $p_n: S(X)\lra
S^n(X)$
is the projection. We have $f_0=\epsilon_C$.

\begin{lem}{} For $n\geq 1$ we have
$$
f_n=\frac{1}{n!}\pi_{S,n}\circ f_1^{\otimes n}\circ \Delta^{(n)}
$$
\end{lem}
\begin{pf}
This follows from the compatibility of $f$ with the comultiplication.
\end{pf}

As a consequence, we get

\begin{prop}{} (\cite{q}, App. B, 4.4) The assignment $f\mapsto f_1$
establishes a bijection between the
set of all unital dg coalgebra maps $f: C\lra S(X)$ and the set of maps
$f_1: C\lra X$
in $\CC (k)$ such that $f_1(1_C)=0$.
\end{prop}

This follows from the previous lemma and the remark that for
$x\in F_iC$ we have $f_n(x)=0$ for $n>i$. $\Box$

\subsubsection{Universal enveloping algebra}
\label{envel} If $\fg$ is a dg Lie algebra,
its {\em universal enveloping algebra $U(\fg)$} is  a dg algebra which
is a quotient of the tensor algebra $T(\fg)$ by the two-sided dg ideal
generated by all elements
$$
xy-(-1)^{|x||y|}yx -[x,y],
$$
$x,y\in \fg=T^1(\fg)$ homogeneous.

The composition $\fg=T^1(\fg)\hra T(\fg)\lra U(\fg)$ is injective;
one identifies $\fg$ with the subcomplex of $U(\fg)$. Evidently,
$U(\fg)$ is generated by $\fg$ as a dg algebra. We shall denote
$F_iU(\fg)\subset U(\fg)$ the subspace spanned by all products of $\leq i$
elements of $\fg$. We have $F_0U(\fg)=k\cdot 1$.

$U(\fg)$ has a canonical structure of a unital cocommutative  dg Hopf
algebra, the comultiplication being defined uniquely by the
requirement that $\fg$ consists of primitive elements. (See ~\cite{q} for
details.) $U(\fg)^+\subset U(\fg)$ will denote the kernel of the counit
({\em augmentation ideal}); one has
a canonical decomposition $U(\fg)\cong k\cdot 1\oplus U(\fg)^+$.

\subsection{Quillen standard complex}
\label{quillen}

\subsubsection{} Let $X\in \CC (k)$. Consider the dg coalgebra $C(X):=S(X[1])$.
Using the d\'{e}calage isomorphisms we identify $C(X)$ with $\oplus_{n\geq 0}
\Lambda^n(X)[n]$.

Let us write down the precise formula for the comultiplication in $C(X)$.
Let $x_i\in X^{\alpha_i},\ i=1,\ldots , n$. For a finite subset
$I=\{ p_1,\ldots , p_i\}\subseteq\{ 1,\ldots , n\};\ p_1<\ldots <p_i$,
set
$$
x_{I}=x_{p_1}\wedge x_{p_2}\wedge\ldots \wedge x_{p_i}\in \Lambda^i(X),
\ x_{\emptyset}=1,
$$
$\bar I:=\{1,\ldots,n\}-I$. Then
\begin{equation}
\label{comult-in-c}
\Delta(x_1\wedge\ldots\wedge x_n)=\sum_I(-1)^{s(I)}x_I\otimes x_{\bar I}
\end{equation}
where the summation is over all subsets $I\subseteq \{1,\ldots,n\}$, and
signs $s(I)$
are defined by the rule
$$
x_1\wedge\ldots\wedge x_n=(-1)^{s(I)}x_I\wedge x_{\bar I}
$$

For future, we shall denote this sign $s(I;\alpha_1,\ldots \alpha_n)$.

Set
$$
C(X)^{pq}=(\Lambda^{-p}(X))^q
$$
(we agree that $\Lambda^p=0$ for $p<0$). It is clear that
$C(X)^n=\oplus_{p+q=n} C(X)^{pq}$. The differential in $C(X)$ has bidegree
$(0,1)$. For future, let us denote it $d_{II}$.
Note that $d_{II}^{-p,*}:\Lambda^p(X)^*\lra \Lambda^p(X)^{*+1}$ is equal to
$(-1)^p$ times the differential on $\Lambda^p(X)$ induced by that on $X$.

\subsubsection{} Now suppose that $\fg$ is a dg Lie algebra. Let us define
maps
$$
d_n: \Lambda^n(\fg)\lra \Lambda^{n-1}(\fg)
$$
by the formula
$$
d_n(x_1\wedge\ldots \wedge x_n)=\sum_{1\leq i<j\leq n}(-1)^{s(\{ i,j\};
\alpha_1,\ldots, \alpha_n)} [x_i,x_j]\wedge x_{\ol{\{ i,j\}}}
$$
where $x_i\in X^{\alpha_i}$.

In particular,

(i) {\em the composition $d_2\circ \pi_{\wedge, 2}$ coincides with the bracket
on $\fg$.}

Set
$$
d_I^{pq}=(d_{-p})^q: C(\fg)^{pq}\lra C(\fg)^{p+1,q};\ d_I=\sum_{p,q} d_I^{pq}
$$

Set $d=d_I+d_{II}$; it is an endomorphism of degree $1$ of the graded
$k$-module $C(\fg)$.

\begin{prop}{} (a) The map $d$ has the following properties.

(ii) $d^2=0$. In particular, $d_I^2=0$ and $d_Id_{II}+d_{II}d_I=0$.

(iii) The comultiplication $\Delta: C(\fg)\lra C(\fg)\otimes C(\fg)$ is
compatible with $d$.

(b) Given $d_{II}$, the differential $d_I$ is uniquely determined by the
properties (i), (ii), and (iii).

(c) Conversely, given $d=d_I+d_{II}$ satisfying (ii) and (iii),
define the bracket $\fg\otimes \fg\lra \fg$ as the composition
$(d_I)^{-2,*}\circ \pi_{\wedge, 2}$. This endows $\fg$ with the structure of a
dg Lie algebra.
\end{prop} $\Box$

So, for a dg Lie algebra $\fg$ we get a unital cocommutative dg coalgebra
$C(\fg)$ which is called
{\em the standard complex} of $\fg$. It was introduced by
Quillen, \cite{q}, App. B, no. 6.

\subsubsection{Maurer-Cartan condition} Let $A$ be a unital cocommutative dg
coalgebra,
and $f: A\lra C(\fg)$ be a map in $\Dgcu$. Let us denote
$f_i: A\lra S^i(\fg[1])$ the composition $p_i\circ f$ where $ p_i$ is
the projection $p_i: C(\fg)\lra S^i(\fg[1])$.

Define the map $[f_1,f_1]: A\lra \fg[2]$ in $\Gr (k)$ as follows:
if $x\in A,\ \Delta (x)=\sum y_i\otimes z_i$, set
$[f_1,f_1](x)=\sum (-1)^{|y_i|} [f_1(y_i),f_1(z_i)]$.

On the other hand, consider
$df_1=d_{\fg}\circ f_1+f_1\circ d_A: A\lra \fg[2]$.

For $x\in A$ we have $f(dx)=f_0(dx)+f_1(dx)+\ldots$. Since
$f(dx)=d\circ f(x)=(d_I+d_{II})\circ f(x)$, we have
$$
f_1(dx)=d_I\circ f_2(x)+d_{II}\circ f_1(x)= d_I\circ f_2(x)-d_{\fg}\circ f_1(x)
$$

On the other hand, if $\Delta (x)=\sum y_i\otimes z_i$, we have by
{}~\ref{coalg-s}
$$
f_2(x)=\frac{1}{2}\sum f_1(y_i)\cdot f_2(z_i)=\frac{1}{2}\sum (-1)^{|y_i|+1}
f_1(y_i)\wedge f_1(z_i)
$$
whence
$$
d_I\circ f_2 =-\frac{1}{2}[f_1,f_1]
$$

Hence we get
\begin{equation}
\label{mc}
df_1+\frac{1}{2}[f_1,f_1]=0
\end{equation}
--- the {\em Maurer-Cartan} equation.

Let us denote $MC(A,\fg)$ the set of all maps $f_1: A\lra \fg[1]$ in $\Gr (k)$
satisfying ~(\ref{mc}) and such that $f_1(1_A)=0$.

\subsubsection{}

\begin{prop}{} (~\cite{q}, App. B, 5.3) The
assignment $f\mapsto f_1$ yields a bijection
$$
\Hom_{\Dgcu}(A,C(\fg))\cong MC(A,\fg)
$$
\end{prop} $\Box$

\subsubsection{}
\label{formula} Suppose that $A$ is a dg Hopf algebra, $a_1,\ldots a_n\in A^0$
primitive elements.

For a subset $I\subset \{1,\ldots,n\}$ set
$a_I=a_{i_1}\cdot a_{i_2}\cdot\ldots\cdot a_{i_s}$ where
$I=\{ i_1,\ldots i_s\}$; $i_1<i_2<\ldots < i_s$; $a_{\emptyset}=1$.

Let us call a {\em $p$-partition} of $\{1,\ldots,n\}$ a sequence
$P=(I_1,\ldots ,I_p)$ of
subsets $I_j\subset \{1,\ldots,n\}$ such that $\{1,\ldots,n\}$ is the
disjoint union $I_1\cup\ldots\cup I_p$. Denote $\CP_p(n)$ the set of all
$p$-partitions.

One computes without difficulty that
$$
\Delta^{(p)}(a_1\cdot\ldots\cdot a_n)=\sum_{P=(I_1,\ldots ,I_p)\in \CP_p(n)}
a_{I_1}\otimes\ldots\otimes a_{I_p}
$$

Suppose we are given a dg coalgebra map $f: A\lra C(\fg)$. It follows form
{}~\ref{coalg-s} that
$$
f_p(a_1\cdot\ldots\cdot a_n)=\frac{1}{p!}
\sum_{P=(I_1,\ldots ,I_p)\in \CP_p(n)}
f_1(a_{I_1})\cdot\ldots\cdot f_1(a_{I_p}) \in S^p(\fg^1)
$$
In particular,
\begin{equation}
f_n(a_1\cdot\ldots\cdot a_n)=f_1(a_1)\cdot\ldots\cdot f_1(a_n)
\end{equation}

\subsection{Connecting morphism}

\subsubsection{Conic dg Lie algebras}
\label{conic} Let $\fg$ be a dg Lie algebra;
$\fh \subset \fg$ a dg Lie ideal. Denote by $i: \fh \lra \fg$ the embedding.

Define a dg Lie algebra $\fX$ as follows. Set $\fX^n=\fh^{n+1}\oplus \fg^n$.
The differential $d:\fX^n\lra \fX^{n+1}$ sends $(h,g)$ to $(-dh,i(h)+dg)$.
So, as a complex, $\fX$ is the usual cone of $i$.

The bracket in $\fX$ is defined as follows. We have $\fX=\fh [1]\oplus \fg$
(as graded modules). The bracket $\fX\otimes \fX\lra \fX$ has components:
$\fg\otimes \fg\lra \fg$ is the bracket in $\fg$;
$\fh [1]\otimes\fg\cong (\fh\otimes\fg)[1]\lra \fh[1]$ and
$\fg\otimes\fh [1]\cong (\fg\otimes\fh)[1]\lra \fh[1]$ are compositions of the
shifting isomorphisms ~(\ref{shift}) and the adjoint action of $\fg$ on $\fh$.
Explicitely:
\begin{equation}
\label{bracket}
[(h,g),(h',g')]=((-1)^a[g,h']+[h,g'],[g,g'])
\end{equation}
for $g\in \fg^a$.

Define maps in $\Gr (k):\ \phi: \fX\lra \fh[1],\ \phi (h,g)=h;\
\theta: \fX\lra \fg,\ \theta (h,g)=g.$ Note that $\phi$ is a map
of complexes. On the other hand, $\theta$ is a morphism
of graded (not dg)  Lie algebras. We have
\begin{equation}
\label{dtheta}
d\theta=-i\circ \phi
\end{equation}

Consequently, $\theta$ induces the map $\Theta: U(\fX)\lra U(\fg)$ of
enveloping algebras which is a morphism of graded (not dg) Hopf algebras.

\subsubsection{Construction of the connecting morphism}
\label{constr} Define the map
$\tilde c_1: T(\fX) \lra \fh[1]$ in $\Gr (k)$ as follows.
Set
$$
\tilde c_1|_{T^0(\fX)}=0,\ \tilde c_1|_{T^1(\fX)}=\phi
$$

Suppose we have defined $\tilde c_1$ on $T^n(\fX)$ for $n\geq 1$.
For $u\in T^n(\fX),\ x\in \fX$ set
\begin{equation}
\label{tc}
\tilde c_1(xu)= (-1)^{|x|}\ad_{\fg}(\theta (x))(\tilde c_1(u))
\end{equation}
This defines $\tilde c_1$ on $T^{n+1}(\fX)$.
Here $\ad_{\fg}$ denotes the adjoint action of $\fg$ on $\fh$:
$\ad_{\fg}(g)(h)=[g,h]$.

Let us denote $\ad_{U(\fg)}$ the induced action of $U(\fg)$ on $\fh$.
{}~(\ref{tc}) implies the important equality:
\begin{equation}
\label{imp}
\tilde c_1(uv)=(-1)^{|u|}\ad_{U(\fg)}(\Theta(u))(\tilde c_1(v))
\end{equation}
for all $u\in T(\fX), v\in T^+(\fX):=\sum_{n>0} T^n(\fX)$.

\subsubsection{}

\begin{thm}{constr-thm} (i) The map $\tilde c_1$ vanishes on the kernel
of the projection $T(\fX)\lra U(\fX)$, and hence it induces the map
\begin{equation}
\label{c1}
c_1: U(\fX)\lra \fh[1]
\end{equation}

(ii) $c_1$ satisfies the Maurer-Cartan equation ~(\ref{mc}).

Consequently, $c_1$ defines the map of unital dg coalgebras
\begin{equation}
\label{conn}
c: U(\fX)\lra C(\fh)
\end{equation}
\end{thm}

We will call $c$ {\bf connecting morphism}.

\begin{pf} (i) We have to prove that
\begin{equation}
\label{toprove}
\tilde c_1(uxyv-(-1)^{ab}uyxv-u[x,u]v)=0
\end{equation}
for all $u,v\in T(\fX),\ x\in \fX^a,\ y\in \fX^b$. If $v\in T^+(\fX)$
this follows from ~(\ref{imp}); so we can suppose $v=1$.

Again, from ~(\ref{imp}) follows that it suffices to prove ~(\ref{toprove})
for $u=1$. This reduces to proving that
$$
(-1)^a[\theta(x),\phi(y)]-(-1)^{(a+1)b}[\theta(y),\phi(x)]=\phi([x,y])
$$
which is equivalent to ~(\ref{bracket}).

(ii) Recall the canonical filtration $F_nU(\fX)$ from ~\ref{envel}.
Let us prove by induction on $n$ that
\begin{equation}
\label{mcprove}
(dc_1+\frac{1}{2}[c_1,c_1])(u)=0
\end{equation}
for all $u\in F_nU(\fX)$. If $n\leq 1$ then both summands in ~(\ref{mcprove})
are $0$.

Suppose we have $x\in \fX^a,\ u\in U(\fX)^b$ such that $\epsilon(u)=0$
where $\epsilon: U(\fX)\lra k$ is the counit. We have
\begin{multline}
(dc_1)(xu)=d(c_1(xu))+c_1(d(xu))=(-1)^ad(\theta(x),c_1(u)])+ c_1(dx\cdot u)+
(-1)^ac_1(x\cdot du)= \\
=(-1)^ad([\theta(x),c_1(u)])+(-1)^{a+1}[\theta(dx),c_1(u)]+[\theta(x),c_1(du)]
\\
\label{eq1}
\end{multline}
Note that
\begin{eqnarray}
[\theta(x),c_1(du)]=[\theta(x),(dc_1)(u)]-[\theta(x),d(c_1(u))]= \nonumber \\
=[\theta(x),(dc_1)(u)]+(-1)^{a+1}d([\theta(x),c_1(u)])+
(-1)^a[d(\theta(x)),c_1(u)] \nonumber
\end{eqnarray}
(we have used that $d([g,h])=[dg,h]+(-1)^a[g,dh]$ for $g\in \fg^a$).
Substituting this in ~(\ref{eq1}) we get
\begin{equation}
\label{dcone}
(dc_1)(xu)=(-1)^a[d(\theta(x))-\theta(dx),c_1(u)]=(-1)^{a+1}[\phi(x),c_1(u)]
\end{equation}
(we have used ~(\ref{dtheta})).

Suppose that $\Delta(u)=u\otimes 1+1\otimes u+\sum_iu'_i\otimes u''_i$.
Let $b,b'_i,b''_i$ be the degrees of $u,u'_i,u''_i$ respectively.
Then
$$
\Delta(xu)=xu\otimes 1+1\otimes xu+x\otimes u+(-1)^{ab}u\otimes x+
\sum_ixu'_i\otimes u''_i+\sum_i(-1)^{ab'_i}u'_i\otimes xu''_i
$$
so that
\begin{multline}
\label{conecone}
[c_1,c_1](xu)=2(-1)^a[\phi(x),c_1(u)]+\sum_i(-1)^{a+b'_i}[c_1(xu'_i),c_1(u''_i)]
+ \\
+\sum_i(-1)^{(a+1)b'_i}[c_1(u'_i),c_1(xu''_i)]=2(-1)^a[\phi(x),c_1(u)]+
\ad_{\fg}(\theta(x))([c_1,c_1](u))
\end{multline}
Adding up ~(\ref{dcone}) and ~(\ref{conecone}) we get
$$
(dc_1+\frac{1}{2}[c_1,c_1])(xu)=\ad_{\fg}(\theta(x))
((dc_1+\frac{1}{2}[c_1,c_1])(u))
$$
Now ~(\ref{mcprove}) follows by induction on $n$.
This proves the theorem.
\end{pf}

\subsubsection{} Let $\fa$ be a graded Lie algebra, $M$ a graded $\fa$-module,
$\phi: \fa\lra M$ a $1$-cocycle of $\fa$ with values in $M$, i.e.
$$
\phi([a,b])=a\phi(b)-(-1)^{|a||b|}b\phi(a)
$$
for $a,b\in \fa$. Then there exists a unique map of graded modules
$$
c_1^+: U(\fa)^+\lra M
$$
such that
(i) $c_1^+|_{\fa}=\phi$; (ii) $c_1^+$ commutes with the action of $\fa$, where
the $\fa$ acts on $U(\fa)^+$ by the left multiplication. Taking the
composition with the projection $U(\fa)\lra U(\fa)^+$ we get a map
$$
c_1: U(\fg)\lra M
$$

The map ~(\ref{c1}) is obtained by applying this remark to
$\fa=\fX,\ M=\fh[1]$,
the action of $\fX$ on $\fh[1]$ being induced from the adjoint action of
$\fg$ on $\fh[1]$ through the graded Lie algebra morphism $\theta:\fX\lra \fg$.

\section{Twisted enveloping algebras and connecting morphism}
\label{twisted}

\subsection{}
\label{dif-smo}
 From now on until the end of the paper we fix
a ground field $k$ of characteristic $0$.

Let $(X,\CO_X)$ be a topological space equipped with a sheaf of
commutative $k$-algebras $\CO_X$. Define the {\em tangent sheaf} $\CT_X$
as the sheaf of $\CO_X$-modules associated with the presheaf
$$
U\mapsto \operatorname{Der}_k(\Gamma(U,\CO_X),\Gamma(U,\CO_X))
$$
(the space of $k$-derivations). $\CT_X$ is a sheaf of $k$-Lie algebras.

We will say that $X$ is {\em differentially smooth} if there exists an open
covering $X=\bigcup U_i$ such that for each $U_i$ the restriction
$\CT_X|_{U_i}$ is a free $\CO_X|_{U_i}$-module admitting a finite basis
of {\em commuting} sections $\dpar_1,\ldots,\dpar_n\in\Gamma(U_i,\CT_X)$.

\subsection{} If $\CF$ a sheaf on $X$, the notation $t\in \CF$ will mean
that $t$ is a local
section of $\CF$. We will use below the straightforward "sheaf" versions
of the definitions and results from Sect.~\ref{envstand}. In particular,
we will use the notion of a {\em  dg $\CO_X$-Lie algebra} (the
bracket is supposed to be $\CO_X$-linear).

If $\CF^{\cdot}$ is a complex of sheaves, $\CZ^i(\CF^{\cdot})$ will denote
sheaves of $i$-cocycles, $\CH^i(\CF^{\cdot})$ cohomology  sheaves
(not to be confused with cohomology {\em spaces} $H^i(X,\CF)$).

If $\fg$ is a dg $\CO_X$-Lie algebra, we will consider its standard
complex $C(\fg)=C_{\CO_X}(\fg)$ (this is a complex of sheaves of dg
$\CO_X$-coalgebras) over $X$) and its canonical filtration $F_iC(\fg)$;
we set $\gr_iC(\fg)=F_iC(\fg)/F_{i-1}C(\fg)$. We will use notations
\begin{equation}
\label{homology}
\CH_i^{Lie}(\fg):=\CH^{-i}(C(\fg));\ \CF_j\CH_i^{Lie}(\fg):=\CH^{-i}(F_jC(\fg))
\end{equation}

\subsection{Dg Lie algebroids} Below we will use some definitions and
constructions from ~\cite{bb}, 1.2 (see also ~\cite{bfm}, 3.2) and their
generalization to a dg-situation.

\begin{defn}{} A {\em dg Lie algebroid} over $X$ is
a sheaf of dg $k$-Lie algebras $\CA$ on $X$ together with a structure
of a left $\CO_X$-module on it and a map
$$
\pi: \CA\lra \CT_X
$$
of dg $k$-Lie algebras and $\CO_X$-modules such that
$$
[a,fb]=f[a,b]+\pi(a)(f)b
$$
for all $a,b\in \CA,\ f\in \CO_X$. We denote $\CA_{(0)}:=\ker \pi$.
It is an $\CO_X$-Lie algebra and a dg Lie ideal in $\CA$.

$\CA$ is called {\em transitive} if $\pi$ is epimorphic (i.e.
$\pi^0: \CA^0\lra \CT_X$ is epimorphic).

A Lie algebroid is a dg Lie algebroid concentrated in degree $0$.
\end{defn}

Dg Lie algebroids over $X$ form a category in an obvious way, with
the final object $(\CT_X, \id_{\CT_X})$.

A dg Lie algebroid with $\pi=0$ is the same as a dg $\CO_X$-Lie algebra.

\subsection{Modules} Let $\CA$ be a dg Lie algebroid over $X$.
A {\em dg $\CA$-module} is a complex of left $\CO_X$ modules $M$ which is
(as a complex of $k$-modules) equipped with the action of the
dg $k$-Lie algebra $\CA$ such that
\begin{equation}
\label{eq1-mod}
f(am)=(fa)m
\end{equation}
\begin{equation}
\label{eq2-mod}
a(fm)=f(am)+\pi(a)(f)m
\end{equation}
for all $f\in \CO_X,\ a\in \CA,\ m\in M$. We will say in this situation
that we have an {\em action} of a dg Lie algebroid $\CA$ on $M$.

\subsection{Twisted enveloping algebras}
\label{twist-env} Let $(\CA, \pi)$ be a dg Lie algebroid
over $X$. Let $U_k(\CA)$ denote the enveloping algebra of $\CA$ considered
as a dg $k$-Lie algebra, $U_k(\CA)^+\subseteq U_k(\CA)$ the augmentation ideal.

Let us consider  $U_k(\CA)^+$ as a sheaf of dg algebras without unit. Denote by
$U_{\CO_X}(\CA)^+$ its quotient by the two-sided dg ideal generated by all
elements
$$
a_1\cdot fa_2-fa_1\cdot a_2-\pi(a_1)(f)a_2,
$$
$a_1,a_2\in \CA,\ f\in \CO_X$.
Set $U_{\CO_X}(\CA)=U_{\CO_X}(\CA)^+\oplus \CO_X$ and define the structure
of a dg algebra with unit by the rule
$$
f\cdot a=fa,\ a\cdot f=fa+\pi(a)(f),\ f\in \CO_X,a\in \CA
$$

We have a canonical algebra map
$$
\CO_X\lra  U_{\CO_X}(\CA)
$$
providing $U_{\CO_X}(\CA)$ with a structure of an $\CO_X$-bimodule, and
a map of left $\CO_X$-modules and dg $k$-Lie algebras (with the structure of
a dg Lie algebra on $U_{\CO_X}(\CA)$ given by the commutator)
$$
i: \CA\lra U_{\CO_X}(\CA)
$$
which is a composition of evident maps $\CA\lra U_k(\CA)^+\lra U_{\CO_X}(\CA)$.
The map $i$ induces an equivalence of the category of dg modules over
$U_{\CO_X}(\CA)$ (as an associative dg algebra) and that of $\CA$-modules.

We define the canonical filtration $F_iU_{\CO_X}(\CA)\subset U_{\CO_X}(\CA)$
as the image of $F_iU_k(\CA)^+\oplus \CO_X$.

\subsection{Coalgebra structure} Let $\CA$ be a dg Lie algebroid
over $X$.  Consider the comultiplication
$\Delta: U_k(\CA)\lra U_k(\CA)\otimes_k U_k(\CA)$. For $x\in U_k(\CA)^+$,
$\Delta(x)=x\otimes 1+1\otimes x+\Delta^+(x)$ where
$\Delta^+(x)\in U_k(\CA)^+\otimes_k U_k(\CA)^+$. Consider the composition
\begin{equation}
\label{comp}
U_k(\CA)^+\overset{\Delta^+}{\lra}U_k(\CA)^+\otimes_k U_k(\CA)^+
\lra U_{\CO_X}(\CA)^+\otimes_{\CO_X} U_{\CO_X}(\CA)^+
\end{equation}
(we consider the tensor square of $U_{\CO_X}(\CA)^+$ as a {\em left}
$\CO_X$-module).

One checks directly that $\ker (U_k(\CA)^+\lra U_{\CO_X}(\CA))$ is contained
in the kernel of ~(\ref{comp}). Hence, ~(\ref{comp}) induces the map
$$
\Delta^+:U_{\CO_X}(\CA)^+\lra U_{\CO_X}(\CA)^+\otimes_{\CO_X} U_{\CO_X}(\CA)^+
$$
Define
$$
\Delta:U_{\CO_X}(\CA)\lra U_{\CO_X}(\CA)\otimes_{\CO_X} U_{\CO_X}(\CA)
$$
as follows: for $f\in \CO_X,\ \Delta(f)=f\otimes 1$; for
$x\in U_{\CO_X}(\CA)^+,\ \Delta(x)=x\otimes 1+1\otimes x+\Delta^+(x)$.
Define the counit $U_{\CO_X}(\CA)\lra \CO_X$ to be the canonical projection.

This defines a structure of a unital cocommutative dg $\CO_X$-coalgebra
on $U_{\CO_X}(\CA)$.

\subsection{Differential operators}

Let  $(X,\CO_X)$ be a differentially smooth $k$-ringed space.

\begin{defn}{}  The sheaf of differential operators
(resp., that of operators of order $\leq n$) on $X$ is
$$
\Diff_X=U_{\CO_X}(\CT_X);\ \Diff^{\leq n}_X=F_nU_{\CO_X}(\CT_X)\ .
$$
\end{defn}

\subsection{Poincare-Birkhoff-Witt condition} The associated graded
algebra
$$
\gr U_{\CO_X}(\CA):=\oplus_{i\geq 0} F_iU_{\CO_X}(\CA)/F_{i-1}U_{\CO_X}(\CA)
$$
is a commutative dg $\CO_X$-algebra. Hence, we have a canonical surjective
morphism
\begin{equation}
\label{pbw-map}
j: S_{\CO_X}{\CA}\lra \gr U_{\CO_X}(\CA)
\end{equation}
We will say that {\em $\CA$ satisfies the Poincare-Birkhoff-Witt (PBW)
condition} if ~(\ref{pbw-map}) is isomorphism.

\subsection{}
\begin{thm}{pbw-algebr}
Let $(X,\CO_X)$ be differentially smooth and let $\CA$ be a transitive dg Lie
algebroid over $X$. Then $\CA$ satisfies PBW.
\end{thm}

{\em Proof} will be done in several steps.

\subsubsection{}
\label{pbw-lie}
\begin{lem}{} Let $\fg$ be a dg $\CO_X$-Lie algebra. Then one has a canonical
isomorphism of unital dg $\CO_X$-coalgebras
$$
e: S_{\CO_X}(\fg)\lra U_{\CO_X}(\fg)
$$
defined by the formula
$$
e(x_1\cdot\ldots\cdot x_n)=\frac{1}{n!}\sum_{\sigma\in \Sigma_n}\pm i(x_1)\cdot
\ldots i(x_n)
$$
where $i:\fg\lra U(\fg)$ is the canonical map, the sign $\pm$ is inserted
according to the Quillen rule.
\end{lem}

\begin{pf} The  argument of the proof of ~\cite{q}, App. B, Thm. 3.2
works in our situation.
\end{pf}

Let us continue the proof of ~\ref{pbw-algebr}. We can forget about
differentials. The question is local on $X$, so we can suppose that
$\CT_X$ is freely generated over $\CO_X$ by $n$ commuting global vector
fields $\dpar_1,\ldots, \dpar_n\in \CT_X(X)$ which can be lifted to
global sections $\gamma_1,\ldots,\gamma_n\in \CA^0(X)$. Consider the map
$$
\mu: U_{\CO_X}(\CA_{(0)})(X)\otimes_k k[\gamma_1,\ldots,\gamma_n]\lra
U_{\CO_X}(\CA)(X)
$$
given by the multiplication. Here $k[\gamma_1,\ldots,\gamma_n]$ denotes the
free left $k$-module on the basis $\gamma_1^{d_1},\ldots,\gamma_n^{d_n}$,
$d_i$ being nonnegative integers.

\subsubsection{}
\begin{lem}{} $\mu$ is an isomorphism of filtered graded $k$-modules.
\end{lem}

Here the filtration on $U_{\CO_X}(\CA_{(0)})(X)$ is the canonical one;
on  $k[\gamma_1,\ldots,\gamma_n]$ the filtration by the total degree,
and the filtration on their tensor product is the tensor product
of filtrations:
$F_i(A\otimes B)=\sum_{p+q=i}\Ima(F_pA\otimes F_qB\lra A\otimes B)$.

\begin{pf} We apply the idea of Serre's proof of the PBW theorem, see
{}~\cite{se}, proof of Thm. 3, p. I, ch.
III.4. First, it is clear that $\mu$ is surjective.

To prove the injectivity,
consider the free left $U_{\CO_X}(\CA_{(0)})(X)$-module $F$ with the basis
$\{ x_M\}$ indexed by finite non-decreasing sequences $M=(i_1,\ldots,i_d)$
with $i_j\in \{1,\ldots,n\}$.
One can introduce on $M$ the structure of an $\CA(X)$-module as follows.

For $M$ as above, call $d$ the length of $M$, and denote it $l(M)$. For
$i\in  \{1,\ldots,n\}$ say that $i\leq M$ if $i\leq i_1$, and denote $iM$ the
concatenation
$(i,i_1,\ldots,i_d)$. Note that $\CA(X)$ is generated as
an abelian group by elements $a\gamma_i,\ a\in \CA_{(0)}(X)$, so it
suffices to define $\gamma_i\cdot ux_M,\ u\in U_{\CO_X}(\CA_{(0)})(X)$.
Let us do this by induction. We set
$$
\gamma_i\cdot ux_M=u\gamma_i\cdot x_M-\ad_{\gamma_i}(u)\cdot x_M
$$
so it suffices to define $\gamma_i\cdot x_M$.
Suppose we have defined $\gamma_i\cdot ux_N$ for all $N:\ l(N)<l(M)$ and
$\gamma_j\cdot ux_N$ for all $j<i,\ l(N)=l(M)$; we suppose that these
elements are linear combinations over $U_{\CO_X}(\CA_{(0)})(X)$ of $x_{N'}$
with $l(N')\leq l(N)+1$. Set
$$
\gamma_i\cdot x_M=\left\{ \begin{array}{ll}
                              x_{iM} & \mbox{if $i\leq M$} \\
                              \gamma_j\cdot (\gamma_i\cdot x_N)+
                                [\gamma_i,\gamma_j]\cdot x_N & \mbox{if $M=jN,\
                                                                      i>j$}
                           \end{array}
                   \right.
$$
the right hand side being defined by induction (note that
$[\gamma_i,\gamma_j]\in \CA_{(0)}(X)$). One checks that this definition is
correct.

Using this, one proves, as in {\em loc. cit.} that  all $x_M$ form the
$U_{\CO_X}(\CA_{(0)})(X)$-basis of $U_{\CO_X}(\CA)(X)$, which implies the
injectivity of $\mu$.
\end{pf}

To finish the proof of ~\ref{pbw-algebr}, it remains to note that
$$
\gr (\mu): S_{\CO_X}(\CA_{(0)})(X)\otimes_k k[\gamma_1\ldots,\gamma_n]\lra
\gr U_{\CO_X}(\CA)(X)
$$
coincides with $j(X)$. Theorem ~\ref{pbw-algebr} is proven. $\Box$

\subsection{}
\label{quasi-iso}
\begin{cor}{}
Let $(X,\CO_X)$ be differentially smooth and $f: \CA\lra \CB$ be a map of
transitive dg Lie algebroids
satisfying one of the two assumprions:

(i) $f$ locally $\CO_X$-homotopy equivalence;

(ii) $f$ is quasiisomorphism, all components $\CA^i,\ \CB^i$ are flat over
$\CO_X$, $\CH^i(\CA)=\CH^i(\CB)=0$ for big $i$.

Then the induced map $U_{\CO_X}(f): U_{\CO_X}(\CA)\lra U_{\CO_X}(\CB)$
is a filtered quasi-isomorphism.
\end{cor}

\begin{pf} Each of our hypotheses implies that
$T_{\CO_X}(f): T_{\CO_X}(\CA)\lra T_{\CO_X}(\CB)$ is a quasi-isomorphism;
hence this is true for $S_{\CO_X}(f)$. Now by ~\ref{pbw-algebr} the same
is true for $U_{\CO_X}(f)$.
\end{pf}

\subsection{Pushout} Let $\CA$ be a dg Lie algebroid over $X$ and
$\fg$ be a dg $\CO_X$-Lie algebra which is also an $\CA$-module.
An {\em $\CA$-morphism} $\psi: \CA_{(0)}\lra \fg$
is a morphism of $\CO_X$-Lie algebras which commutes with the $\CA$-action,
where the action of $\CA$ on $\CA_{(0)}$ is the adjoint one, and such that
$$
a\cdot x=[\psi(a),x]
$$
for all $a\in \CA_{(0)},\ x\in \fg$.

Given such a morphism, the dg Lie algebroid $\CA_{\psi}$ is defined as follows.
Consider the dg Lie algebra semi-direct product $\CA \semid \fg$ (so, the
bracket is
$[(a,x),(b,y)]=([a,b],-(-1)^{|b||x|}b\cdot x+a\cdot y+[x,y]),\ a,b\in \CA,\
x,y\in fg$). By definition, $\CA_{\psi}$ is the quotient of $\CA \semid \fg$
by the dg Lie ideal $\CA_{(0)}\hra \CA \semid \fg,\ a\mapsto (a, -\psi(a))$.
The map $\pi_{\CA_{\psi}}:\CA_{\psi}\lra \CT_X$ maps $(a,x)$ to $\pi_{\CA}(a)$.

If $\CA$ is transitive then so is $\CA_{\psi}$, and $\CA_{\psi (0)}=\fg$.

\subsection{Boundary morphism}

\subsubsection{Conic Lie algebroids}
Let $(\CA,\pi)$ be a dg Lie algebroid over $X$, $\fh$ a dg Lie algebra and a
left $\CO_X$-module, $i: \fh\hra \CA$ an embedding of dg Lie algebras and
of $\CO_X$-modules such that $i(\fh)$ is a dg Lie ideal in $\CA$, and
$\pi\circ i=0$. This implies that $\fh$ is a dg $\CO_X$-Lie algebra.

Let us consider the complex $\fA=\Cone (i)$. According to ~\ref{conic},
$\fA$ has a canonical structure of a dg $k$-Lie algebra. Together with the
evident structure of an $\CO_X$-module and
$\pi_{\fA}:\fA\overset{\theta}{\lra} \CA\overset{\pi}{\lra} \CT_X$,
$\fA$ becomes a dg Lie algebroid.

Let us consider $\fA$ as a dg $k$-Lie algebra, and apply to it (the
sheaf version of) the construction ~\ref{constr}. We get a map of
 sheaves of graded $k$-modules
\begin{equation}
\label{conn-prep}
c_{1/k}^+: U_k(\fA)^+\lra \fh[1]
\end{equation}
satisfying the Maurer-Cartan equation.

\subsubsection{}
\begin{thm}{}
\label{conn-thm} The map $c_{1/k}^+$ factors through the
canonical map $U_k(\fA)^+\lra U_{\CO_X}(\fA)^+$ and hence induces the map
$$
c_1^+:U_{\CO_X}(\fA)^+\lra \fh[1].
$$
Taking its composition with the projection
$U_{\CO_X}(\fA)\lra U_{\CO_X}(\fA)^+$, we get
$$
c_1:U_{\CO_X}(\fA)\lra \fh[1].
$$
This map is a morphism of left graded $\CO_X$-modules and it satisfies
the Maurer-Cartan equation.
Consequently, it induces the map of filtered dg $\CO_X$-coalgebras
\begin{equation}
c:U_{\CO_X}(\fA)\lra C_{\CO_X}(\fh)
\end{equation}
\end{thm}

The map $c$ will be called {\bf boundary morphism} associated to $\fA$.

\begin{pf} We have only to check that
\begin{equation}
\label{check}
ux(fy)v-u(fx)yv-u(\pi_{\fA}(x)(f)y)v\in \ker(c_{1/k}^+)
\end{equation}
for all $x,y\in \fA,\ u,v\in U_k(\fA),\ f\in \CO_X$.
Note that $\fh$ is an $\CA$-module with respect to the adjoint action
(NB! this is not true for $\CA$: the Axiom ~(\ref{eq1-mod}) does not hold).
It follows that the action of $U_k(\CA)$ on $\fh[1]$ factors through
$U_{\CO_X}(\CA)$, hence ~(\ref{check}) holds true for $v\in U_k(\fA)^+$;
so we can suppose $v=1$; obviously, it is enough to check ~(\ref{check})
for $u=1$.In that case it is a direct check (note that the maps $\phi$ and
$\theta$ are $\CO_X$-linear):
$$
c_{1/k}^+(x(fy)-(fx)y-\pi_{\fA}(x)(f)y)=[\theta(x),f\phi(y)]-
[f\theta(x),\phi(y)]-\pi_{\CA}(\theta(x))(f)\phi(y)=0
$$
\end{pf}

\subsection{Connecting morphisms}
\label{abstr-ksmaps}

Here $X$ is supposed to be differentially smooth.

\subsubsection{} Let $\CA$ be a transitive dg Lie algebroid
over $X$. Set $\fh:=\CA_{(0)}$, so we have an exact sequence
\begin{equation}
\label{fund-abstr}
0\lra \fh\overset{i}{\lra} \CA\overset{\pi}{\lra} \CT_X\lra 0\ .
\end{equation}
It induces the map
\begin{equation}
\label{ks1-abstr}
\kappa^1: \CT_X\lra \CH^1(\fh)\ .
\end{equation}

Set $\fA=\Cone (i)$. By Thm.~\ref{quasi-iso} $\pi$ induces a
filtered quasi-isomorphism $U_{\CO_X}(\fA)\lra U_{\CO_X}(\CT_X)=\Diff_X$,
whence isomorphisms
$$
\CH^0(U_{\CO_X}(\fA))\cong \Diff_X;\ \CH^0(F_nU_{\CO_X}(\fA))\cong
\Diff_X^{\leq n}\ .
$$
On the other hand, Thm.~\ref{conn-thm} gives the filtered map
$$
c: U_{\CO_X}(\fA)\lra C_{\CO_X}(\fh)\ .
$$
By taking $\CH^0(c)$, we get maps
\begin{equation}
\label{abstr-ks}
\kappa:\Diff_X\lra \CH_0^{Lie}(\fh)
\end{equation}
as well as
\begin{equation}
\label{abstr-ks-i}
\kappa^{\leq n}: \Diff^{\leq n}_X\lra \CF_n\CH_0^{Lie}(\fh)
\end{equation}
which are called {\bf connecting morphisms}.

\subsubsection{} We have by definition
$$
\gr_nC_{\CO_X}(\fh)=S^n_{\CO_X}(\fh[1])
$$
hence the maps
\begin{equation}
\label{proj-lie}
\CF_n\CH_0^{Lie}(\fh)\lra \CH_0(S^n_{\CO_X}(\fh[1]))\ .
\end{equation}

On the other hand, the embeddings $S^n(\fh^1)\subset S^n(\fh[1])$ induce
embedding of cocycles $S^n\CZ^1(\fh)\hra \CZ^0(S^n(\fh[1]))$ which pass to
cohomology and give the maps
\begin{equation}
\label{emb}
S^n_{\CO_X}(\CH^1(\fh))\lra \CH^0(S^n(\fh[1]))\ .
\end{equation}

\subsubsection{}
\label{main-thm}
\begin{thm}{} The connecting morphisms ~(\ref{abstr-ks}) and
{}~(\ref{abstr-ks-i}) have the following properties.

(i) The squares
$$\begin{array}{ccc}
\Diff^{\leq n-1}_X & \overset{\kappa^{\leq n-1}}{\lra} &
\CF_{n-1}\CH_0^{Lie}(\fh) \\
\downarrow & \; & \downarrow \\
\Diff^{\leq n}_X & \overset{\kappa^{\leq n}}{\lra} & \CF_{n}\CH_0^{Lie}(\fh) \\
\end{array}$$
commute. We have
$$
\kappa=\lim_{\ra}\kappa^{\leq n}\ .
$$

(ii) The squares
$$\begin{array}{ccccc}
\Diff^{\leq n}_X &\; & \overset{\kappa^{\leq n}}{\lra} &\;
& \CF_{n}\CH_0^{Lie}(\fh) \\
\downarrow & \; & \; &\; & \downarrow \\
S^n_{\CO_X}(\CT_X) & \overset{(-1)^nS^n(\kappa^1)}{\lra} & S^n_{\CO_X}
(\CH^1(\fh)) & \overset{(\ref{emb})}{\lra} & \CH^0(S^n(\fh[1])) \\
\end{array}$$
commute. Here the left vertical arrow is the symbol map,
and the right one is ~(\ref{proj-lie}).
\end{thm}

The property (i) is obvious. (ii) will be proven in the next Subsection.

\subsection{Explicit formulas}

\subsubsection{}
\label{expl1} In the previous assumptions, suppose we  have local
sections $\dpar_1,\ldots ,\dpar_n\in \CT_X$.
Let us pick $0$-cocycles in
$\fA$ lifting them: $a_p=(-\alpha_p,\gamma_p)\in \fA^0,\ \gamma_p\in
\CA^0,\alpha_p\in\fh^1,\ \pi(\gamma_p)=\dpar_p;\
d_{\CA}(\gamma_p)=i(\alpha_p)$.

Consider the map $c_{1/k}^+$, ~(\ref{conn-prep}). For $I$ as in ~\ref{formula}
define $a_I$ as there, and set
$\dpar_I=\dpar_{i_1}\cdot\ldots\cdot\dpar_{i_s}$,
$$
\alpha(\dpar_I)=c_{1/k}^+(a_I)=\ad (\gamma_{i_1})\circ\ldots\circ
\ad(\gamma_{i_{s-1}})(-\alpha_{i_s})
$$
Define elements (we use the notations of ~\ref{formula})
$$
\kappa^{(p)}_n=\frac{1}{p!}
\sum_{P=(I_1,\ldots,I_p)\in\CP_p(n)}\alpha(\dpar_{I_1})\cdots
\alpha(\dpar_{I_p})\in S^p_{\CO_X}(\fh^1),
$$
$p=1,\ldots,n$. In particular,
\begin{equation}
\label{kappa-nn}
\kappa^{(n)}_n=(-1)^n\alpha_1\cdots
\alpha_n.
\end{equation}

\subsubsection{}
\label{expl2}
\begin{prop}{} The class $\kappa^{\leq n}(\dpar_1\cdots\dpar_n)$
is represented by the cocycle
$$
(0,\kappa^{(1)}_n,\ldots,\kappa^{(n)}_n)\in\oplus_{p=0}^nS^p_{\CO_X}(\fh^1)
\subset F_nC_{\CO_X}(\fh)^0.
$$
\end{prop}

\begin{pf} This follows from ~\ref{formula} applied to $c_{1/k}$ (note
that $\fA$ is a Hopf algebra).
\end{pf}

As a corollary, we get the claim (ii) of ~\ref{main-thm} which follows
from ~(\ref{kappa-nn}).

\subsubsection{Schur polynomials} Let us define polynomials
$P_n(\alpha_1,\alpha_2,\ldots)$ by means of the generating function
$$
\exp (\sum_{p=1}^{\infty}\alpha_p\frac{t^p}{p!})=
\sum_{n=0}^{\infty}P_n(\alpha_1,\ldots)\frac{t^n}{n!}
$$
It is easy to see that $P_n(\alpha_1,\ldots)=P_n(\alpha_1,\ldots,\alpha_n)$
and
\begin{equation}
\label{schur}
P_n(\alpha_1,\ldots,\alpha_n)=\sum_{(n_1,\ldots,n_s):\sum_j jn_j=n}
\frac{n!}{(1!)^{n_1}\cdot\ldots\cdot (s!)^{n_s}
n_1!\cdot\ldots\cdot n_s!} \alpha_1^{n_1}\cdot\ldots\cdot\alpha_s^{n_s}
\end{equation}

The first polynomials are: $P_0=1,\ P_1=\alpha_1,\ P_2=\frac{\alpha_1^2}{2}+
\alpha_2$.

\subsubsection{} Now suppose we are given a local section $\dpar\in \CT_X$
together with a lifting $a=(-\alpha,\gamma)\in \fA^0$ as in ~\ref{expl1}.
For $i\geq 1$ set
$$
\alpha_i=(\ad (\gamma))^{i-1}(-\alpha)\in \fh^1
$$

\subsubsection{}
\label{explicit}
\begin{prop}{} The class $\kappa^{\leq n}(\dpar^n)$ is represented by the
cocycle
$$
P_n(\alpha_1,\ldots,\alpha_n)\in F_nS_{\CO_X}(\fh^1)\subseteq
F_nC_{\CO_X}(\fh)^0.
$$
\end{prop}

\begin{pf} From ~\ref{expl2} follows that $\kappa^{\leq n}(\dpar^n)$ can be
represented by the polynomial $Q(\alpha_1,\ldots,\alpha_n)$ where
the coefficient of $Q$ at $\alpha_1^{n_1}\cdots\alpha_s^{n_s}$
is equal to the number of partitions of the set $\{1,\ldots,n\}$ containing
$n_1$ of $1$-element subsets, $n_2$ of $2$-element subsets, ...,
$n_s$ of $s$-element subsets. From ~(\ref{schur}) follows that
$Q=P_n$.
\end{pf}

\subsubsection{} Summing up the expressions ~\ref{explicit} over $n$, we can
rewrite our formulas as
$$
\kappa(\exp(t\dpar))=\exp (\frac{\exp(\ad(\gamma))-\Id}{\ad(\gamma)}(\alpha)).
$$
This was pointed out to us by I.T.Leong. Cf. Deligne's formula
{}~\cite{gmd}, p. 51, (1-1).



\section{Differential calculus on formal schemes}
\label{formal}

The aim of this Section is to extend (a part of) the classical Grothendieck's
language of differential calculus, \cite{ega} IV, \S\S 16, 17, to formal
schemes.

\subsection{} Recall that we have fixed a ground field $k$ of
characteristic $0$.
We will work in the category $\Fsch$ of separated locally noetherian formal
schemes over $k$. All necessary definitions and facts about them are contained
in \cite{ega} I, \S10. Objects of $\Fsch$ will be called simply
{\em formal schemes}.
Inside $\Fsch$, we will consider a full subcategory $\Sch$ of separated
locally noetherian schemes of over $k$ whose objects will be called simply
(usual) {\em schemes}.

By definition, a formal scheme is a topological space $\fX$ equipped with
a sheaf of topological rings $\CO_{\fX}$. $\fX$ is a union of
affine formal schemes $\Spf(A)$ where $A$ is a noetherian $k$-algebra
complete in the $I$-adic topology for some ideal $I\subset A$. As a topological
space, $\Spf(A)=\Spec(A/I)$, and $\Gamma(\Spf(A),\CO_{\Spf(A)})=A$.

There exists a sheaf of ideals $\fI\subset\CO_{\fX}$ such that for sufficiently
small affine $U\subset \fX$ the ideals $\Gamma(U,\fI)^n,\ n\geq 1$, form
a base for the topology of $\Gamma(U,\CO_{\fX})$. Such a sheaf is called
the {\em  ideal of definition} of $\fX$. All ringed spaces
$$
X_n:=(\fX,\CO_{\fX}/\fI^{n+1})
$$
are (usual) schemes, and we have
\begin{equation}
\label{indlim}
\fX=\dirlim\ X_n
\end{equation}
cf. {\em loc.cit.}, 10.11. Among ideals of definition there exists a unique
maximal one.

Let $f:\fX\lra \fS$ be a morphism of formal schemes, $\fK\subseteq\CO_{\fS}$
an ideal of definition of $\fS$ and $\fI$ the maximal ideal of definition
of $\fX$. Then $f^*(\fK)\subseteq\fI$, so $f$ induces maps of schemes
$f_n:X_n\lra S_n$ such that
\begin{equation}
\label{adic}
f=\dirlim\ f_n
\end{equation}
cf. {\em loc.cit.}, 10.6.10.

\subsubsection{} The category $\Fsch$ has fibered products, {\em loc.cit.},
10.7.

\subsection{}
\label{aff}  Let $\fX$ be affine, $\fX=\Spf(A)$. We have a canonical
functor
$$
\Delta: \Mod(A)\lra \Mod(\CO_{\fX}),\ M\mapsto M^{\Delta}
$$
from the category of $A$-modules to the category of sheaves
$\CO_{\fX}$-modules.

If $A$ is noetherian then $\Delta$ establishes an equivalence between
the category of $A$-modules of finite type and that of coherent $\CO_{\fX}$-
modules, {\em loc.cit.}, 10.10.2.

\subsection{} Let $\fX$ be a formal scheme, $\CM$ a coherent
$\CO_{\fX}$-module.
If $\fX$ is represented as in ~(\ref{indlim}) then
\begin{equation}
\label{indlim-m}
\CM\cong \invlim\ M_n
\end{equation}
for a suitable inverse system of coherent $X_n$-modules $M_n$, {\em loc.cit.},
10.11.3.

We will consider $\CM$ as a sheaf of topological $\CO_{\fX}$ modules equipped
with the topology defined in {\em loc.cit.}, 10.11.6. Note that
for every affine open $U\subseteq\fX$ the module $\Gamma(U,\CM)$ is
complete.

If $\CN$ is another coherent $\CO_{\fX}$-module, we have isomorphisms
$$
\CM\otimes_{\CO_{\fX}}\CN\cong\invlim\ (M_n\otimes_{\CO_{X_n}}N_n)
$$
and
$$
\CHom_{\CO_{\fX}}(\CM,\CN)\cong \invlim\ \CHom_{\CO_{X_n}}(M_n,N_n),
$$
cf. {\em loc.cit.}, 10.11.7.

\subsection{Jets} Let $f: \fX\lra \fS$ be a morphism of formal schemes,
$f=\dirlim\ f_i$ a representation as in ~(\ref{adic}).

Let us consider the diagonal
$$
\Delta_f: \fX\hra \fX\times_{\fS}\fX
$$
We have
$$
\fX\times_{\fS}\fX=\dirlim\ X_i\times_{S_i}X_i
$$
and $\Delta_f=\dirlim\ \Delta_{f,i}$ where
$$
\Delta_{f,i}: X_i\hra X_i\times_{S_i}X_i
$$
are diagonal mappings. If
$\CI_i\subseteq\CO_{X_i\times_{S_i}X_i}$ is the ideal of $\Delta_{f,i}$
then
$$
\fI=\invlim\CI_i\subseteq\invlim\ \CO_{X_i\times_{S_i}X_i}=
\CO_{\fX\times_{\fS}\fX}
$$
is the ideal of $\Delta_f$.

Let $p_j:\fX\times_{\fS}\fX\lra \fX,\ j=1,2$ be projections.

For any integer $n\geq 0$ set
$$
\fX^{(n)}_f=(\Delta_f(\fX),\CO_{\fX\times_{\fS}\fX}/\CI^{n+1})
$$
--- it is a closed formal subscheme of $\fX\times_{\fS}\fX$.
Consider canonical projections $p_i^{(n)}:\fX^{(n)}\lra \fX,\ i=1,2$.

\subsubsection{}
\label{defn-jets}
\begin{defn}{} (Cf. \cite{ega} IV 16.3.1, 16.7.1.) We define a sheaf of
rings over $\fX$ which is called the {\em sheaf of $n$-jets} as
$$
\CP^n_f=\CP^n_{\fX/\fS}:=p_{1*}(\CO_{\fX\times_{\fS}\fX}/\fI^{n+1})=
p_{1*}^{(n)}p^{(n)*}_2(\CO_{\fX})
$$
We have two canonical morphisms of sheaves of topological rings
\begin{equation}
\label{first}
\CO_{\fX}\lra \CP_f^n,\ x\mapsto x\otimes 1
\end{equation}
and
\begin{equation}
\label{second}
d^n_{\fX/\fS}:\CO_{\fX}\lra \CP_f^n,\ x\mapsto 1\otimes x
\end{equation}
by means of which one introduces a structure of left (resp., right)
$\CO_{\fX}$-module on $\CP^n_f$.

More generally, for a coherent sheaf $\CM$ over $\fX$ set
$$
\CP^n_{\fX/\fS}(\CM):=p_{1*}^{(n)}p^{(n)*}_2(\CM)
$$
\end{defn}

We have
$$
\CP^n_{\fX/\fS}(\CM)=\CP^n_{\fX/\fS}\otimes_{\CO_{\fX}}\CM
$$
where the right $\CO_{\fX}$-module structure on $\CP^n_{\fX/\fS}$
is used on the right-hand side, cf. {\em loc.cit.}, 16.7.2.1.
The $\CO_{\fX}$-bimodule structure on $\CP^n_f$ induces an
 $\CO_{\fX}$-bimodule structure on $\CP^n_{\fX/\fS}(\CM)$.

Following {\em loc.cit.}, if we do not specify the structure of an
$\CO_{\fX}$-module on $\CP^n_{\fX/\fS}(\CM)$, we mean that
of a {\em left} module.

\subsubsection{} If $\CM=\invlim\ M_i$ as in ~(\ref{indlim-m}) then
\begin{equation}
\label{invlim-p}
\CP^n_{\fX/\fS}(\CM)=\invlim\ \CP^n_{X_i/S_i}(M_i)\ .
\end{equation}

\subsubsection{} Consider affine open formal subschemes
$U=\Spf(B)\subseteq\fX,\ V=\Spf(A)\subseteq \fS$
such that $f(U)\subseteq V$, so that $A,B$ are adic rings, and $B$ is a
topological $A$-algebra. Let $I=\ker(B\otimes_AB\lra B)$ be the kernel
of the multiplication, and
$$
P_{B/A}=(B\otimes_AB)/I^{n+1}
$$
considered as a $B$-module by means of a map $x\mapsto x\otimes 1$.
It is naturally a topological $B$-module, and we can consider its completion,
$\hP_{B/A}$. We have
\begin{equation}
\CP^n_U\cong (\hP_{B/A})^{\Delta}
\end{equation}

\subsubsection{}
\begin{claim}{} $d^n_{\fX/\fS}(\CO_{\fX})$ topologically generates
$\CO_{\fX}$-module $\CP^n_f$.
\end{claim}

(cf. \cite{ega} IV, 16.3.8) $\Box$

\subsubsection{} We have evident projections
\begin{equation}
\label{cofilt}
\CP^n_f\lra\CP^{n-1}_f\ .
\end{equation}

Consider a projection
\begin{equation}
\label{proj}
\CP^n_f\lra \CP^0_f=\CO_{\fX}
\end{equation}
and let $\bCP^n_f$ denote its kernel. A map
$$
d_1: \CO_{\fX}\lra \CP^n_f,\ x\mapsto x\otimes 1,
$$
is left inverse to ~(\ref{proj}), so we have a splitting
$$
\CP^n_f\cong\bCP^n_f\oplus\CO_{\fX}.
$$

\subsection{Differentials} We define
$$
\Omega^1_f=\Omega^1_{\fX/\fS}:=\bCP^1_f\ .
$$
If $\fX=\Spf(B),\ \fS=\Spf(B)$ then we will use the notation $\Omega^1_{B/A}$
for $\Gamma(\fX,\Omega^1_{\fX/\fS})$.
We have natural isomorphisms
$$
\Omega^1_f\cong p_{1*}(\CI/\CI^2)\cong p_{2*}(\CI/\CI^2)\ .
$$
We have a canonical continuous $\CO_{\fS}$-linear map
\begin{equation}
\label{difl}
d:\CO_{\fX}\lra \Omega^1_{\fX/\fS}
\end{equation}
induced by the map $x\mapsto x\otimes 1-1\otimes x$.

\subsubsection{}
 From ~(\ref{invlim-p}) follows that we have a natural isomorphism
\begin{equation}
\Omega^1_{\fX/\fS}\cong\invlim\ \Omega^1_{X_i/S_i}
\end{equation}

\subsection{Example}
\label{power} Suppose $\fS=\Spf(A)$, and $\fX=\Spf(A_n)$ where
$A_n=A\{ T_1,\ldots,T_n\}$ is the completion of the polynomial
ring $A[T_1,\ldots,T_n]$ in the $J$-adic topology $J$, being an ideal
generated by an ideal
of definition of $A$ and $T_1,\ldots, T_n$. Then $\Omega^1_{A_n/A}$ is a
free $A_n$-module with the basis $dT_1,\ldots,dT_n$.

\subsection{Derivations}  Let $\CM$ be a coherent $\CO_{\fX}$-module. We
define {\em the sheaf of derivations}
$$
\Der_{\fS}(\CO_{\fX},\CM):=\CHom_{\CO_{\fX}}(\Omega^1_{\fX/\fS},\CM).
$$
The map ~(\ref{difl}) induces a canonical embedding
$$
\Der_{\fS}(\CO_{\fX},\CM)\hra\CHom_{\CO_{\fS}}(\CO_{\fX},\CM)
$$
which identifies $\Der_{\fS}(\CO_{\fX},\CM)$ with the sheaf of local
$\CO_{\fS}$-homomorphisms $\alpha:\CO_{\fX}\lra\CM$ such that
$$
\alpha(xy)=x\alpha(y)+y\alpha(x).
$$

We set
$$
\CT_{\fX/\fS}:=\Der_{\fS}(\CO_{\fX},\CO_{\fX})
$$
and call this sheaf the {\em tangent sheaf of $\fX$ with respect to $\fS$},
cf. {\em loc.cit.}, 16.5.7. It is a sheaf of $\CO_{\fS}$-Lie algebras.

\subsection{Differential operators} Let $\CM,\ \CN$ be coherent
$\CO_{\fX}$-modules.
We define {\em the sheaf of differential operators of order $\leq n$ from
$\CM$ to $\CN$}:
$$
\Diff^{\leq n}_{\fX/\fS}(\CM,\CN)=\CHom_{\CO_{\fX}}(\CP^n_{\fX/\fS}(\CM),\CN).
$$
The structure of $\CO_{\fX}$-bimodule on $\CP^n({\CM})$ induces the structure
of $\CO_{\fX}$-bimodule on $\Diff^{\leq n}_{\fX/\fS}(\CM,\CN)$.

The projections ~(\ref{cofilt}) induce embeddings
$$
\Diff^{\leq n-1}_{\fX/\fS}(\CM,\CN)\hra\Diff^{\leq n}_{\fX/\fS}(\CM,\CN)\ .
$$
We set
$$
\Diff_{\fX/\fS}(\CM,\CN)=\dirlim\ \Diff^{\leq n}_{\fX/\fS}(\CM,\CN)\ .
$$

We denote
$$
\Diff^{\leq n}_{\fX/\fS}:=\Diff^{\leq n}_{\fS}(\CO_{\fX},\CO_{\fX});\
\Diff_{\fX/\fS}:=\Diff_{\CO_{\fS}}(\CO_X,\CM)\ .
$$

As in the case of schemes, we have compositions
$$
\Diff^{\leq n}_{\fX/\fS}\otimes_{\CO_{\fS}}\Diff^{\leq m}_{\fX/\fS}\lra
\Diff^{\leq n+m}_{\fX/\fS}
$$
(where in the tensor product the 1st (resp., 2nd) factor has a structure
of a right (resp., left) $\CO_{\fX}$-module).

\subsection{} Note that we have canonically
$$
\Diff^{\leq 1}_{\fX/\fS}\cong\CO_{\fX}\oplus\CT_{\fX/\fS}\ .
$$
Thus, the multiplication induces a canonical map
\begin{equation}
\label{twist-dif}
U_{\CO_{\fX}}(\CT_{\fX/\fS})\lra \Diff_{\fX/\fS}
\end{equation}
--- cf. ~\ref{twist-env}.

\subsection{}
\label{functor} Suppose we have a commutative square of formal schemes
$$\begin{array}{ccc}
\fX & \overset{u}{\lra} & \fX' \\
\downarrow & \; & \downarrow \\
\fS & \lra & \fS'
\end{array}$$
It induces natural maps
$$
\nu_n: u^*\CP^n_{\fX'/\fS'}\lra \CP^n_{\fX/\fS}
$$
If the square is cartesian, these maps are isomorphisms, cf. {\em loc.cit.},
16.4.5.

\subsection{}
\label{point} Let $\fX$ be a formal scheme, $x:\Spec(k)\lra \fX$ a point.
Then one has a canonical isomorphism of $k(x)$-algebras
$$
(\CP^n_{\fX/k})_x\otimes_{\CO_x}k(x)\cong \CO_x/\fm_x^{n+1}
$$
--- cf. {\em loc.cit.}, 16.4.12.

\subsection{Formally smooth morphisms} Let us call a morphism $f:\fX\lra \fS$
of formal schemes {\em formally smooth} if for every commutative
square
$$\begin{array}{ccc}
Y & \lra & \fX \\
\downarrow & \; & \downarrow \\
T & \lra & \fS
\end{array}$$
where $Y=\Spec(B)\lra T=\Spec(A)$ is a closed embedding of affine
(usual) schemes
corresponding to an epimorphism $A\lra B$ whose kernel $I$ satisfies $I^2=0$,
there exists a lifting $T\lra\fX$.

\subsubsection{} Suppose that $f=\dirlim\ f_n$, where all
$f_n:X_n\lra S_n$ are smooth morphisms of (usual) schemes. Then $f$ is
formally smooth.

This follows from ~\cite{ega}  0$_{\mbox{IV}}$, 19.4.1.

\subsection{} Let $\fX\overset{f}{\lra}\fY\overset{g}{\lra}\fS$ be
morphisms of formal schemes.

By ~\ref{functor} they induce maps
$f^*\Omega^1_{\fY/\fS}\lra\Omega^1_{\fX/\fS}$ and
$\Omega^1_{\fX/\fS}\lra\Omega^1_{\fX/\fY}$.
They in turn induce the maps in the sequence
\begin{equation}
\label{seq}
0\lra\CT_{\fX/\fY}\overset{g'}{\lra}\CT_{\fX/\fS}
\overset{f'}{\lra} f^*\CT_{\fY/\fS}\ .
\end{equation}

\subsubsection{}
\label{seq-thm}
\begin{thm}{} The sequence ~(\ref{seq}) is exact. If $f$ is formally
smooth then $f'$ is surjective.
\end{thm}

\begin{pf} Follows from ~\cite{ega} 0$_{\mbox{IV}}$, 20.7.18.
\end{pf}

\subsection{Smooth morphisms}
\label{smooth}

\begin{defn}{} A morphism $f:\fX\lra\fS$ is called {\em smooth} if it is
formally smooth and locally of finite type (\cite{ega} I, 10.13).
\end{defn}

\subsubsection{}
\begin{lem}{} If $f$ is smooth then $\Omega^1_f$ is locally free of finite
type.
\end{lem}
\begin{pf}
Everything is reduced to the case of affine formal schemes,
$\fX=\Spf(B),\ \fS=\Spf(A)$, and $B$ is a noetherian topological $A$-algebra
topologically of finite type over $A$. First, $\Omega^1_{B/A}$ is of finite
type by ~\ref{power} and ~\ref{seq-thm} since $B$ is a quotient of some
$A\{ T_1,\ldots, T_n\}$.

Let us prove that $\Omega^1_{B/A}$ is projective $B$-module. Suppose we have
a diagram
$$\begin{array}{ccc}
\; & \; & \Omega^1_{B/A}\\
\; & \; & \downarrow\\
M & \overset{\psi}{\lra}& N\\
\end{array}$$
with epimorhic $\psi$ and finitely generated $M$ . Let $I$ be an ideal of
definition in $B$.
By \\ \cite{ega}~0$_{\mbox{IV}}$, 20.4.9,
we can find liftings $\phi_n:\Omega^1_{B/A}\lra M/I^nM$ for each $n$.

We claim that we can choose $\phi_n$'s in such a way that
for all $n$ the composition $\Omega^1_{B/A}\overset{\phi_{n+1}}{\lra}
M/I^{n+1}M\lra M/I^nM$ is equal to $\phi_n$. Indeed, we do it step by step,
using again the lifting property and the fact that canonical maps
$$
M/I^{n+1}M\lra M/I^nM\times_{N/I^nN}N/I^{n+1}N
$$
are surjective.

Since $B$ is complete, $M=\invlim\ M/I^nM$, hence there exists a lifting
$\phi:\Omega^1_{B/A}\lra M$.

Now, since $\Omega^1_{B/A}$ is of finite type, we can find an epimorphic
map $F\lra \Omega^1_{B/A}$ from a finitely generated free $B$-module $F$;
using the lifting property proved above, we conclude that $\Omega^1_{B/A}$
is a direct summand of $F$.
\end{pf}

\subsection{Differentially smooth morphisms}

\begin{defn}{} A morphism $f:\fX\lra\fS$ of formal schemes is called
{\em differentially smooth} if $\Omega^1_f$ is locally free of finite rank.
\end{defn}

For an arbitrary $f$, define the graded ring
$$
\CGr.(\CP_f)=\oplus_{n=0}^{\infty}\CGr_n(\CP_f);\
\CGr_n(\CP_f)=\fI^n/\fI^{n+1},
$$
$\fI$ being as in ~\ref{defn-jets}. Evidently $\CGr_1(\CP_f)=\Omega^1_f$ and
we have the canonical surjective morphism from the symmetric algebra
\begin{equation}
\label{symm}
S^._{\CO_{\fX}}(\Omega^1_f)\lra\CGr_.(\CP_f)\ .
\end{equation}

\subsubsection{}
\begin{lem}{} If $f$ is differentially smooth then ~(\ref{symm}) is
isomorphism.
\end{lem}

\begin{pf} The same as in ~\cite{ega} IV 16.12.2 (cf. {\em loc.cit.},
16.10).
\end{pf}

\subsubsection{}
\label{basis}
\begin{thm}{} Let $f$ be differentially smooth, $U\subseteq\fX$ an open,
$\{ t_i\}_{i\in I}\in\Gamma(U,\CO_{\fX})$ a set of sections such that
$\{ dt_i\}_{i\in I}$ is the basis of $\Omega^1_{U/\fS}$. Let
$\{\dpar_i\}\in\Gamma(U,\CT_{\fX/\fS})$ be the dual basis. Then

(i) all $\dpar_i$ commute with each other;

(ii) $\Diff^{\leq n}_{\fX/\fS}$ is freely generated by all monomials on
$\dpar_i$ of degrees $\leq n$.
\end{thm}

\begin{pf} The same as in {\em loc.cit.}, 16.11.2.
\end{pf}

Thus, if $\fS=\Spec(k)$, if $f$ is differentially smooth then
$\fX$ is differentially smooth in the sense of ~\ref{dif-smo}.

\subsection{}
\begin{cor}{} If $f$ is differentially smooth then the map ~(\ref{twist-dif})
is isomorphism.
\end{cor}
\begin{pf}
Follows immediately from \Thm{basis} and \Thm{pbw-algebr}: the
map~(\ref{twist-dif}) preserves filtrations and induces an isomorphism on
the associated graded rings.
\end{pf}

Hence, if $\fS=\Spec(k)$, our sheaf $\Diff_{\fX/\fS}$ coincides with
$\Diff_{\fX}$ defined in Section ~\ref{twisted}.

\subsection{Examples of differentially smooth morphisms}
\subsubsection{} Smooth morphisms. This follows from {\em loc.cit.}, 16.10.2.

\subsubsection{} The structure morphism
$\Spf(A\{ T_1,\ldots,T_n\})\lra\Spf(A)$ (see ~\ref{power}).

\section{Homotopy Lie algebras and direct image functor}
\label{direct}

In this Section we develop a formalism of homotopy Lie algebras which
is sufficiently good for our needs.

In~\ref{holie} we define the category $\Holie$ of homotopy Lie algebras
and a functor from it to the filtered derived category given by the Quillen
standard complex (see~\ref{quillen}). In order to define higher direct
images for $\Holie$ in~\ref{dihla}, we provide in~\ref{s-derham} a construction
of Thom-Sullivan functor from cosimplicial modules to complexes, and some
homotopical properties of it in~\ref{cob}. The proof of the properties of
the Thom-Sullivan functor is given in Section~\ref{thoms}.
The main result of this Section is~\Thm{main-5}.

\subsection{Homotopy Lie algebras}
\label{holie}

\subsubsection{}
\label{cat}
Let $X$ be a formal scheme. We define $\Dglie(X)$
(resp., $\Dglie^{qc}(X),\ \Dglie^{c}(X)$) as  a category
of dg $\CO_X$-Lie algebras $\fg$ such that ---

(i) all components $\fg^i$ are $\CO_X$-flat
(and quasicoherent or coherent over $\CO_X$ respectively);

(ii) one has $\CH^i(\fg)=0$ for sufficiently big $i$.

A morphism $f: \fg\lra \fh$ in this category
is a map of complexes of $\CO_X$-modules compatible with
brackets. Let us call $f$ {\em a quasi-isomorphism} if it induces an
isomorphism of all cohomology sheaves $\CH^i(f):\CH^i(\fg)\lra \CH^i(\fh)$.

By definition, the category $\Holie^{(*)}(X)$, is
the localization of $\Dglie^{(*)}(X)$ with respect to the
class of all quasi-isomorphisms. Objects of $\Holie(X)$ are called
{\em homotopy Lie algebras}. So, each dg $\CO_X$-Lie algebra defines
a homotopy Lie algebra.

\subsubsection{} If $f: X\lra Y$ is a flat morphism, the inverse
image functor for $\CO_X$-modules induces the functor
\begin{equation}
\label{inv-im-dglie}
f^{*Lie}:\Dglie^{(*)}(Y)\lra \Dglie^{(*)}(X)
\end{equation}
It takes quasi-isomorphisms to quasi-isomorphisms, and hence induces
the functor (to be denoted by the same letter)
\begin{equation}
\label{inv-im-holie}
f^{*Lie}:\Holie^{(*)}(Y)\lra \Holie^{(*)}(X)
\end{equation}

\subsubsection{Filtered derived categories} (For more details, see ~\cite{i},
ch. V, no.1 where the case of finite filtrations is considered.)
Let $\CC(X)$ denote the category of complexes of $\CO_X$-modules,
and $\CCF(X)$ the category
whose objects are complexes of $\CO_X$-modules $A$ together with a
filtration by $\CO_X$-subcomplexes $\ldots\subseteq F_iA\subseteq
F_{i+1}A\subseteq\ldots,\ i\in \Bbb Z$ such that $F_iA=0$ for sufficiently
small $i$ and
$A=\cup_i F_iA$; the morphisms being morphisms of complexes
compatible with filtrations.

We set $\gr_i(A)=F_iA/F_{i-1}A$. For $a,b\in \Bbb Z$, let $\CCF_{[a,b]}(X)$
be the full subcategory consisting of complexes
with $F_{a-1}A=0,F_{b}A=A$.

A morphism $f:A\lra B$ in $\CCF(X)$ is called {\em a filtered
quasi-isomorphism} if the induced maps $\gr_i(f): \gr_i(A)\lra \gr_i(B)$ are
quasi-isomorphisms.    From our assumptions on filtrations follows that
a filtered quasi-isomorphism is a quasi-isomorphism.

We denote by $\CD(X)$ the localization of $\CC(X)$ with respect
to quasi-isomorphisms, and $\CDF(X)$,  $\CDF_{[a,b]}(X)$ the localization
of $\CCF(X)$ (resp., $\CCF_{[a,b]}(X)$) with respect to  filtered
quasi-isomorphisms.

\subsubsection{}
\begin{lem}{} Suppose that $f:\fg\lra \fh$ is a quasi-isomorphism in
$\Dglie(X)$. Then the induced morphism $C(f): C(\fg)\lra C(\fh)$ is a filtered
quasi-isomorphism.
\end{lem}

\begin{pf} It suffices to prove that all
$\gr_i(f): S^i(\fg[1])\lra S^i(\fh[1])$ are quasi-isomorphisms.
Repeated application of the K\"{u}nneth spectral sequence ~\cite{g},
Ch.~I, 5.5.1,
shows that $f^{\otimes i}:\fg^{\otimes i}\lra \fh^{\otimes i}$ are
quasi-isomorphisms (we use the assumptions ~\ref{cat} (i)-(ii)).
After passing to $\Sigma_n$-invariants, we get the desired claim.
\end{pf}

\subsubsection{} It follows that the functors $\fg\mapsto C(\fg)$,
$\fg\mapsto F_nC(\fg)$ induce functors between homotopy categories
$$
C: \Holie(X)\lra \CDF(X)
$$
$$
F_nC:\Holie(X)\lra \CDF_{[0,n]}(X)
$$
For $\fg\in \Holie(X)$ we define homology sheaves
$$
\CH^{Lie}_n(\fg)=\CH^{-n}(C(\fg));\ \CF_m\CH_n^{Lie}(\fg)=H^{-n}(F_mC(\fg))
$$
--- these are sheaves of $\CO_X$-modules. If $X=\Spec (k)$, we denote them
$H^{Lie}_n(\fg),\ F_mH_n^{Lie}(\fg)$ respectively --- these are
$k$-vector spaces.

\subsubsection{}
\label{sec-spectr-seq} Let $\fg\in \Holie(X)$. We have a spectral sequence
\begin{equation}
\label{spectr-seq}
E_1^{pq}=\CH^q(\Lambda^{-p}_{\CO_X}(\fg)) \Lra \CH^{Lie}_{-p-q}(\fg)
\end{equation}
If all sheaves $\CH^i(\fg)$ are $\CO_X$-flat then by K\"{u}nneth formula,
we have
\begin{equation}
\label{kunneth}
\CH^q(\Lambda^{n}_{\CO_X}(\fg))\cong
\left(\sum_{q_1+\ldots+q_n=q}\CH^{q_1}(\fg)\otimes_{\CO_X}\ldots\otimes_{\CO_X}
\CH^{q_n}(\fg)\right)^{\Sigma_n,-}
\end{equation}
where $(\cdot)^{\Sigma_n,-}$ denotes the subspace of anti-invariants of
the symmetric group $\Sigma_n$.

\subsection{Thom-Sullivan complex}
\label{s-derham}

\subsubsection{} In this subsection $\Delta$ will denote the category
of totally ordered finite sets $[n]=\{ 0,\ldots, n\}$,
$n\geq 0,$
and non-decreasing maps. For any category $C$, we denote $\Delta^0C,\
\Delta C$ the categories of simplicial and cosimplicial objects
in $C$ respectively. For $A\in \Delta^0 C$ (resp., $B\in \Delta C$) we
denote $\alpha^*: A_m\lra A_n$ (resp., $\alpha_*: A^n\lra A^m$) the
map in $C$ corresponding to a morphism $\alpha: [n]\lra[m]$
in $\Delta$.

$\Ens$ will denote the category of sets, $\Delta[n]\in \Delta^0\Ens$ the
standard $n$-simplex.

\subsubsection{} (For details, see \cite{bug}).
Denote by $R_n$ the commutative $k$-algebra
{}~{$k[t_0,\ldots,t_n]/(\sum_{i=0}^nt_i-1)$} ($t_i$ being independent
variables).
Set $\Bbb A[n]=\Spec R_n$.
Together with the standard face and degeneracy maps, algebras $R_n,\ n\geq 0,$
form a simplicial algebra $R$.

Let $\Omega_n$ denote the algebraic de Rham complex of $R_n$ over $k$,
$$
\Omega_n=\Gamma(\Bbb A[n],\Omega_{\Bbb A[n]/k})
$$
It is a commutative dg $k$-algebra which may be identified with
$$
R_n[dt_0,\ldots,dt_n]/(\sum dt_i),
$$
with $\deg(dt_i)=1$ and the differential given by the formula $d(t_i)=dt_i$.
The algebras $\Omega_n,\ n\geq 0$ together with coface and codegeneracy
maps induced from $R$, form a simplicial  commutative dg algebra
$\Omega=\{\Omega_n\}$. So, for a fixed $p,\ \Omega^p$ is a
simplicial vector space, and $\Omega_p$ is a complex.

\subsubsection{}
\label{invtensor} (Cf. ~\cite{hlha}, 3.1.3) Let $C$ be a small category.
Let us denote $\Mor(C)$
the category whose objects are morphisms $f: x\lra y$ in $C$,
a map $f\lra g$ is given by a commutative diagram
$$\begin{array}{ccc}
\cdot & \overset{f}{\lra} & \cdot \\
\uparrow & \; & \downarrow \\
\cdot & \overset{g}{\lra} & \cdot \\
\end{array}$$
The composition is defined in the evident way. For $f$ as above set $s(f)=x,\
t(f)=y$.

Let $A$ be a commutative $k$-algebra, $X: C^0\lra \Mod(k),\
Y:C\lra \Mod(A)$ two functors ($C^0$ denotes the opposite
category).  Denote by
$$
X\otimes Y: \Mor(C)\lra \Mod(A)
$$
the functor given by $X\otimes Y(f)=X(s(f))\otimes_kY(t(f))$ and defined
on morphisms in the evident way. Set
$$
X\invtimes Y=\invlim (X\otimes Y)\in \Mod(A)
$$
In other words, $X\invtimes Y$ is an $A$-submodule of
$\prod_{c\in Ob(C)}X(c)\otimes Y(c)$ consisting of all $\{a(c)\}_{c\in Ob(C)},$
$ a(c)\in X(c)\otimes Y(c)$ such that for every
$f:b\lra c\in \Mor(C)$
$$
(f^*\otimes 1)(a(c))=(1\otimes f_*)(a(b)).
$$

\subsubsection{Thom-Sullivan complex} (Cf. ~\cite{hlha}, \S 4) Applying
the previous construction to $C=\Delta,\ X=\Omega^p\in \Delta^0\Mod(k),
\ Y\in\Delta\Mod(A)$, we get $A$-modules
$$
\Omega^p(Y)=\Omega^p\invtimes Y
$$
When $p$ varies, they form a complex $\Omega(Y)\in \CC(A)$ which is called
{\em Thom-Sullivan complex} of $Y$.

If $Y$ is a constant cosimplicial object, then
\begin{equation}
\label{const}
\Omega(Y)\cong Y
\end{equation}
canonically.

For example, if $Y$ is the cosimplicial space of $k$-valued cochains of a
simplicial set $X$, $\Omega(Y)$ is the Thom-Sullivan complex of
$X$ described in ~\cite{bug}.

If $Y$ is a {\em complex} in $\Delta\Mod(A)$ or, what is the same,
a cosimplicial complex of $A$-modules, then, applying the previous
construction componentwise we get a {\em bicomplex} of $A$-modules;
we will denote the corresponding simple complex again $\Omega(Y)$.
This way we get a functor
\begin{equation}
\label{thom}
\Omega: \Delta\Mod(A)\lra \CC(A).
\end{equation}

\subsubsection{}
\begin{lem}{} The functor $\Omega$ is exact.
\end{lem}

For a proof, see \ref{exact-6}.

\subsubsection{Normalization} For $Y\in \Delta\Mod(A)$ denote by
$N(Y)\in \CC(A)$ its
normalization, i.e. set $N(Y)^i\subset Y^i$ to be the intersection
of kernels of all codegeneracies $Y^i\lra Y^{i-1}$, the differential
$N(Y)^i\lra N(Y)^{i+1}$ being the alternating sum of the cofaces. ($N(Y)^i=0$
for $i<0$). We say that $Y$ is {\em finite dimensional} if $N(Y)^i=0$
for $i>>0$.
This way we get a functor
$$
N:\Delta\Mod(A)\lra \CC(A)
$$

For each $n\geq 0$ denote by $\CZ^{\cdot}_n=C^{\cdot}(\Delta[n],k)$ the
complex of
normalized $k$-valued cochains of the standard simplex. When $n$ varies,
the complexes $\CZ^{\cdot}_n$ form a simplicial object
$\CZ=\CZ^{\cdot}_{\cdot}\in \Delta^0\CC(k)$.

Given $Y\in \Delta\Mod(A)$, we can apply the construction
{}~\ref{invtensor} to each $\CZ^n_{\cdot}$ and $Y$, and obtain
$\CZ\invtimes Y\in \CC(A)$. It follows from the definitions that
we have a natural isomorphism
$$
\CZ\invtimes Y\cong N(Y)
$$
(cf. ~\cite{hlha}, 2.4).

\subsubsection{} For each $n$ we have the {\em integration}
map
$$
\int: \Omega_n\lra C^*(\Delta[n],k)
$$
in $\CC(k)$, ~\cite{bug},\S 2. Taken together, they give rise to the morphism
$$
\int: \Omega\lra \CZ
$$
in $\Delta^0C(k)$. It induces natural maps
\begin{equation}
\label{int-y}
\int_Y:\Omega(Y)\lra N(Y)
\end{equation}
for every $Y\in \Delta\Mod(A)$.

\subsubsection{"De Rham theorem"}
\label{derham-thm}
\begin{lem}{} For every $Y\in \Mod(A)$ the map $\int_Y$ is a
quasi-isomorphism.
\end{lem}
\begin{pf}
 See ~\cite{hlha}, 4.4.1.
\end{pf}

\subsubsection{Base change} Let $A'$ be a commutative $A$-algebra,
$Y\in \Delta\Mod(A)$, whence $Y\otimes_AA'\in \Delta\Mod(A')$. We have an
evident base change morphism
\begin{equation}
\label{bchange-map}
\Omega(Y)\otimes_AA'\lra \Omega(Y\otimes_AA').
\end{equation}

\subsubsection{}
\label{base}
\begin{lem}{} If $Y$ is finite dimensional then ~(\ref{bchange-map}) is
isomorphism.
\end{lem}

\subsubsection{}
\label{flat}
\begin{lem}{} Let $A$ be noetherian and $Y\in \Delta\Mod(A)$;
suppose that all $Y^n$ are flat over $A$.
Then for every $p,\ \Omega^p(Y)$ is flat over $A$.
\end{lem}

The proof of ~\ref{base} and ~\ref{flat} will be given in the next Section,
see ~\ref{base-6},~\ref{flat-6}.

\subsubsection{}
\label{cosimpl-lie} Suppose that $B$ is a commutative dg algebra and
$\fg$ is a dg Lie algebra. Then $B\otimes \fg$ is naturally
a dg Lie algebra, the bracket being defined as
$$
[a\otimes x,b\otimes y]=(-1)^{|x||b|} ab\otimes [x,y].
$$

Let $\fg$ be a cosimplicial dg $A$-Lie algebra. Then all
$\Omega_n\otimes \fg^n$ are dg $A$-Lie algebras; hence their inverse limit
$\Omega(\fg)$ is. This way we get a functor $\Omega$ from the
category of cosimplicial dg $A$-Lie algebras to dg $A$-Lie algebras.

\subsubsection{} Let $X$ be a (formal) scheme. We can sheafify Thom-Sullivan
construction. Denote by $\Mod(\CO_X)$ the category of sheaves
of $\CO_X$-modules. For $\CF\in \Delta\Mod(\CO_X)$ we get
$N(\CF), \Omega(\CF)\in \CC(X)$. We call $\CF$ {\em finite dimensional}
if $N(\CF)^i=0$
for $i>>0$. We denote $\Delta(\Mod(\CO_X))^f\subseteq \Delta(\Mod(\CO_X))$ the
full subcategory consisting of finite dimensional cosimplicial sheaves.

By~\ref{qc} if $\CF$ is finite dimensional, and all $\CF^i$ are quasicoherent
then all $\Omega^i(\CF)$ are quasicoherent.

By~\ref{flat} if all $\CF^i$ are $\CO_X$-flat then all $\Omega^i(\CF)$
are $\CO_X$-flat.

If $\fg$ is a cosimplicial dg $\CO_X$-Lie algebra, then, applying
{}~\ref{cosimpl-lie} we get a dg $\CO_X$-Lie algebra $\Omega(\fg)$.

\subsubsection{}
\label{derham-dir-im} Let $f:X\lra Y$ be a map of schemes,
$\CF\in\Delta\Mod(\CO_X)$. Then we have an evident equality
$$
f_*\Omega(\CF)=\Omega(f_*\CF)
$$

\subsection{Cosimplicial homotopies}
\label{cob}
Let $\CA$ be a category with finite products.
\subsubsection{}
\begin{defn}{path}
 Define the {\em path
functor} $X\mapsto X^I$ from $\Delta{\CA}$ to itself as follows:
$$ (X^I)^n=\prod_{s:[n]\ra[1]}X^n,$$
the map $f_*:(X^I)^m\ra(X^I)^n$ for any $f:[m]\ra[n]$ being defined
by the formula
$$ f_*(\{x_s\})_t=f_*(x_{tf}),\ t:[n]\ra[1].$$
\end{defn}
The path functor is endowed with natural transformations
$$i_X:X\ra X^I,\ \pr_{i,X}:X^I\ra X (i=0,1)$$
given by the formulas
$$ i_X(x)_s=x;\ \pr_{i,X}(\{x_s\})=x_{s_i}$$
where $s_i:[n]\ra[1]$ denotes the constant map with value $i$.

\Defn{path} is a special case of constructions given in~\cite{q2}, ch.~2,
see Prop. 2 for the dual statement. In particular, one can define in the same
way a functor $X\mapsto X^S$ for any finite simplicial set $S$, this
construction is functorial on $S\in\Delta^0\Ens$ and the natural
transformations $i_X,\pr_{i,X}$ are induced by the corresponding maps
$\pi:I\ra *\text{ and }\iota_i:*\ra I (i=0,1)$ in $\Delta^0\Ens$.

\subsubsection{}
\begin{defn}{strict-homo}
Let $X,Y\in\Delta{\CA}$. Maps $f_i:X\ra Y \ (i=0,1)$ are said to be
{\em strictly homotopic} if there exists a map $F:X\ra Y^I$ such that
$f_i=\pr_i\circ F$.
\end{defn}

\subsubsection{}
\begin{exa}{=}
Let $X\in\Delta{\CA}$. The maps $\id,\ i\circ\pr_0:X^I\ra X^I$ are
strictly homotopic.

In fact, the maps $\id_I$ and $\iota_0\circ\pi$ are strictly homotopic
in $\Delta^0\Ens$.
\end{exa}

\subsubsection{} Let now $\CA$ be additive. Recall that
for $X\in\Delta{\CA}$ the total  complex $\Tot(X)\in \CC({\CA})$
is defined by the
properties
$$ \Tot(X)^n=X^n,$$
$$d=\sum (-1)^i\delta^i:\Tot(X)^n\ra\Tot(X)^{n+1}.$$
This defines  a functor $\Tot:\Delta{\CA}\ra \CC({\CA})$.

\begin{lem}{homo-homo}
A strict homotopy $H:X\ra Y^I$ between $f$ and $g$ induces a (chain)
homotopy $h:\Tot(X)\ra\Tot(Y)[-1]$ between $\Tot(f)$ and $\Tot(g)$.
\end{lem}
\begin{pf}
For $x\in X^n$ define
$$ h(x)=\sum_{i=0}^{n-1}(-1)^i\sigma^i(y_{\alpha_i})$$
where $y_{\alpha}$ are defined by $H(x)=\{y_{\alpha}\}$ and
\begin{equation}
\alpha_i(t)=\begin{cases} 0& \text{ if } t\leq i\\
    1& \text{ if } t>i.
            \end{cases}
\end{equation}
A direct calculation shows that $h$ is the chain homotopy we need.
\end{pf}

\subsubsection{} We apply the above constructions to homotopy Lie algebras.
Fix a commutative ring $A\supseteq{\Bbb Q}$ and put ${\CA}=\Dglie(A)$.

\begin{cor}{glueh}
Let $X,Y\in\Delta\Dglie(A)$. Let $f,g:X\ra Y$ be strictly homotopic. Then
the maps $\Omega(f), \Omega(g):\Omega(X)\ra\Omega(Y)$ induce equal maps
in the homotopy category $\Holie(A)$.
\end{cor}
\begin{pf}
It suffices to check that $\Omega(\pr_{0,Y})=\Omega(\pr_{1,Y})$ in the homotopy
category. Since the both maps split $\Omega(i_Y)$, it sufficies to check that
the latter one is a quasi-isomorphism. For this we can substitute the
functor $\Omega$ with $\Tot$ (since they are naturally quasi-isomorphic).
We can also substitute the category $\Dglie(A)$ with the category ${\CC}(A)$
since the forgetful functor $\#:\Dglie(A)\ra{\CC}(A)$ commutes with direct
products. Then~\Exa{=} and~\Lem{homo-homo} prove even more that we actually
need.
\end{pf}

\subsection{Direct image of homotopy Lie algebras}
\label{dihla}
\subsubsection{\v{C}ech resolutions} Let $X$ be a topological space, $\CF$
an abelian sheaf over $X$, $\CU=\{ U_i\}_{i\in I}$ an open covering of $X$.
For each $n\geq 0$ set
$$
\CHC^n(\CU,\CF)=\prod_{(i_0,\ldots,i_n)\in I^{n+1}}j_{i_0\ldots i_n*}
j_{i_0\ldots i_n}^*\CF
$$
where $j_{i_0\ldots i_n}:U_{i_0}\cap\ldots\cap U_{i_n}\hra X$.
Together with the standard cofaces and codegeneracies, the sheaves
$\CHC^n(\CU,\CF),\ n\geq 0,$ form a cosimplicial sheaf $\CHC(\CU,\CF)$.
It is finite dimensional if the covering $\CU$ is finite.

We have an embedding $\CF\lra \CHC^0(\CU,\CF)$; it induces the augmentation
map
\begin{equation}
\label{augm}
\CF\lra\CHC(\CU,\CF)
\end{equation}
where $\CF$ is considered as a constant cosimplicial sheaf. The induced
map
\begin{equation}
\label{augm-n}
\CF\lra N(\CHC(\CU,\CF))
\end{equation}
is a quasi-isomorphism, \cite{g}, ch.~II, 5.2.1.

\subsubsection{} Suppose that

--- either $X$ is a scheme and $\CF$ is a quasicoherent sheaf of
$\CO_X$-modules,

--- or $X$ is a formal scheme and $\CF$ is a coherent sheaf.

Choose an affine covering $\CU$. Then the complex
$\Gamma(X,N(\CHC(\CU,\CF)))$ represents $R\Gamma(X,\CF)$,  ~\cite{h},
III, 4.5; ~\ref{aff}.

If $f:X\lra Y$ is a morphism of (formal) schemes then $f_*N(\CHC(\CU,\CF))$
represents $Rf_*(\CF)$.

\subsubsection{} Applying to ~(\ref{augm}) the functor $\Omega$, and using
{}~(\ref{const}), we get a canonical map of complexes
\begin{equation}
\label{augm-omega}
\CF\lra \Omega(\CHC(\CU,\CF))
\end{equation}
which is a quasi-isomorphism by the above and ~\ref{derham-thm}.

\subsubsection{} Suppose we are in one of the following situations:

{\bf Case 1.}  $f: X\lra Y$ is a flat morphism of schemes,
$\fg\in \Dglie^{qc}(X)$.

{\bf Case 2.} $f: X\lra Y$ is a flat morphism of formal schemes,
$\fg\in \Dglie^{c}(X)$.

Choose an open affine covering $\CU$ of $X$,
and consider the cosimplicial complex of $\CO_Y$-modules
$f_*\CHC(\CU,\fg)$; it has an evident structure of a cosimplicial
dg $\CO_Y$-Lie algebra. Applying the Thom-Sullivan functor,
we get a dg $\CO_Y$-Lie algebra $\Omega(f_*\CHC(\CU,\fg))$.
Note that by ~\ref{derham-dir-im} this is the same as
$f_*\Omega(\CHC(\CU,\fg))$. Let us denote it $f_{*,\CU}^{Lie}(\fg)$.

\subsubsection{} Let now  $\CU=\{U_i\}_{i\in I}$ and
$\CV=\{V_j\}_{j\in J}$ be two open coverings of $X$  so that cosimplicial
dg Lie algebras $\CHC(\CU,\fg),\ \CHC(\CV,\fg)$ are defined.
Let maps $f,g:I\ra J$
satisfy the conditions $U_i\subseteq V_{f(i)},\ U_i\subseteq V_{g(i)}.$

\begin{lem}{to5.}
The maps from $\CHC(\CV,\fg)$ to $\CHC(\CU,\fg)$ induced by $f$ and $g$
are strictly homotopic.
\end{lem}
\begin{pf} Immediate.
\end{pf}

\begin{cor}{} In the notations above the maps from $f_{*,\CV}^{Lie}(\fg)$
to $f_{*,\CU}^{Lie}(\fg)$ induced by $f$ and by $g$, coincide in $\Holie(Y)$.
\end{cor}
\begin{pf}
Compare~\ref{to5.} with~\ref{glueh}.
\end{pf}

\subsubsection{} If $\fg=f^{*\Lie}\fh$ for some $\fh\in\Dglie(Y)$ then the
augmentation ~(\ref{augm-omega})
$$
f^{*\Lie}\fh\lra \Omega(\CHC(\CU,f^{*\Lie}\fh))
$$
and the adjunction map $\fh\lra f_*f^*\fh$ induce the map
of dg $\CO_X$-Lie algebras
\begin{equation}
\label{adj}
\fh\lra f_{*,\CU}^{Lie}(f^{*\Lie}\fh).
\end{equation}

On the other hand, for any $\fg\in\Holie(X)$ a map
\begin{equation}
f^{*\Lie}f_*^{\Lie}(\fg)\lra\fg
\end{equation}
in $\Holie(X)$ is defined as the composition
\begin{equation}
f^{*\Lie}f_*^{\Lie}(\fg)=f^{*\Lie}f_*\Omega(\CHC(\CU,\fg))\ra
\Omega(\CHC(\CU,\fg))\ra\fg
\end{equation}
the last map being inverse to the composition of the
quasi-isomorphisms~(\ref{augm-n}) and~(\ref{int-y}).

Putting together the above considerations, we get the following

\subsubsection{}
\begin{thm}{main-5} (i) In Case 1 (resp., Case 2) $f_{*,\CU}^{Lie}(\fg)$
belongs to
$\Dglie^{qc}(Y)$ (resp., $\Dglie(Y)$).

(ii) The class of $f_{*,\CU}^{Lie}(\fg)$ in $\Holie^{qc}(Y)$ (resp.,
$\Holie(Y)$) does not depend, up to a unique isomorphism, on $\CU$.

(iii) The functor $f_{*,\CU}^{Lie}$ takes quasi-isomorphisms to
quasi-isomorphisms; thus it induces the functor
$$
f_*^{Lie}:\Holie^{qc}(X)\lra \Holie^{qc}(Y)
$$
in Case 1, and
$$
f_*^{Lie}:\Holie^{c}(X)\lra \Holie(Y)
$$
in Case 2 respectively, such that the square
$$\begin{array}{ccc}
\Holie^{(*)}(X) & \overset{f_*^{Lie}}{\lra} & \Holie^{(*')}(Y) \\
\downarrow & \; & \downarrow \\
\CD(X)     & \overset{f_*}{\lra} & \CD(Y) \\
\end{array}$$
the vertical arrows being forgetful functors, 2-commutes in both cases.
This means that there is a natural isomorphism between the two functors
from $\Holie^{(*)}(X)$ to $\CD(Y)$.

(iv) In Case 1, maps ~(\ref{adj}) induce the natural transformation
$\Id_{\Holie(Y)}\lra f_*^{Lie}f^{*Lie}$ which makes the functor $f_*^{Lie}$
right adjoint to $f^{*Lie}$. $\Box$
\end{thm}

\subsubsection{} In case $Y=\Spec(k),\ f:X\lra Y$ the structure
morphism, we will denote $f_*^{Lie}(\fg)$ also by $\Gamma^{Lie}(X,\fg)$.

\section{Thom-Sullivan functor}
\label{thoms}

In this Section we will compute more explicitely the Thom-Sullivan
functor and prove some fundamental properties of it.

\subsection{} We keep the assumptions and notations from ~\ref{s-derham}.
In particular, $k$ is a base field of characteristic zero and
$A$ denotes a commutative $k$-algebra.

Expressing $t_0$ as $t_0=1-\sum_{i=1}^nt_i$ we identify $k$-algebras
$R_n$ with $k[t_1,\ldots,t_n]$ and commutative
dg $k$-algebras $\Omega_n$ with $R_n[dt_1,\ldots,dt_n]$.

The standard simplicial morphisms of dg $k$-algebras
$$
d_i:\Omega_n\lra\Omega_{n-1},\ s_i:\Omega_n\lra
\Omega_{n+1},
$$
$i=0,\ldots, n$ are defined by the formulas

\begin{equation}
d_i(t_j)=
  \begin{cases}
     t_j & \text{ if } j<i,\\
     0   & \text{ if } j=i,\\
     t_{j-1}& \text{ if } j>i
  \end{cases}
\label{d_i}
\end{equation}
for $i>0$ and
\begin{equation}
d_0(t_j)=
  \begin{cases}
     t_{j-1} & \text{ if } j>1,\\
     1-\sum t_i& \text{ if } j=1.
  \end{cases}
\label{d_0}
\end{equation}

\begin{equation}
s_i(t_j)=
  \begin{cases}
     t_j & \text{ if } j<i,\\
     t_i+t_{i+1}  & \text{ if } j=i,\\
     t_{j+1}& \text{ if } j>i
  \end{cases}
\label{s_i}
\end{equation}
for $i>0$ and
\begin{equation}
s_0(t_j)=t_{j+1}.
\label{s_0}
\end{equation}

These maps satisfy standard simplicial identities.

\subsection{} Let $X\in\Delta\Mod(A)$. So  $X$ is a set of $A$-modules
$X^p,\ p\geq 0$, together with maps
$$
\delta^i:X^{p-1}\lra X^{p},\ \sigma^i:X^{p+1}\lra X^{p},
$$
$i=0,\ldots,p$, satisfying the standard cosimplicial identities.

By definition, $\Omega(X)$ is a complex of $A$-modules
$$
0\lra\Omega^0(X)\overset{d}{\lra}\ldots\overset{d}{\lra}\Omega^n(X)
\overset{d}{\lra}\ldots
$$
where $\Omega^n(X)$ is the space of all collections
$\{x_p\in\Omega_p^n\otimes X^p\}_{p\geq 0}$ satisfying the following conditions
\begin{equation}
(1\otimes\delta^i)(x_p)=(d_i\otimes 1)(x_{p+1})
\label{d-eq}
\end{equation}

\begin{equation}
(s_i\otimes 1)(x_p)=(1\otimes\sigma^i)(x_{p+1})
\label{s-eq}
\end{equation}
for all $p,i$ --- see the diagram below.


\begin{center}
  \begin{picture}(14,6)
    \put(5,0){\makebox(4,2){$\Omega^n_{p+1}\otimes X^{p+1}$}}
    \put(5,4){\makebox(4,2){$\Omega^n_p\otimes X^{p}$}}
    \put(0,2){\makebox(4,2){$\Omega^n_{p+1}\otimes X^{p}$}}
    \put(10,2){\makebox(4,2){$\Omega^n_p\otimes X^{p+1}$}}

    \put(5.5,4.5){\vector(-2,-1){2}}
    \put(5.5,1.5){\vector(-2,1){2}}
    \put(8.5,1.5){\vector(2,1){2}}
    \put(8.5,4.5){\vector(2,-1){2}}

   \put(3.5,4){\makebox(1,0.5){$\scriptsize s_i\otimes 1$}}
   \put(3.5,1.5){\makebox(1,0.5){$\scriptsize 1\otimes\sigma^i$}}
   \put(9.5,1.5){\makebox(1,0.5){$\scriptsize d_i\otimes 1$}}
   \put(9.5,4){\makebox(1,0.5){$\scriptsize 1\otimes\delta^i$}}
  \end{picture}
\end{center}


Until the end of this Section, let us fix $n\geq 0$. Our aim will be
an explicit computation of $\Omega^n(X)$.

\subsection{} In what follows ${\Bbb N}=\{0,1,\ldots\}$.

Let $p\in\Bbb N$. For $a\in{\Bbb N}^p,\alpha\in\{0,1\}^p$ denote
$$
\omega_{a,\alpha}=t_1^{a_1}\cdots t_p^{a_p}dt_1^{\alpha_1}\wedge
\cdots\wedge dt_p^{\alpha_p}\in\Omega_p.
$$
We have $\deg\omega_{a,\alpha}=\sum\alpha_i$.

Let $I^p\subset {\Bbb N}^p\times\{0,1\}^p$ consist of all pairs
$(a,\alpha)$ such that $\sum\alpha_i=n$.

Evidently, forms $\omega_{a,\alpha},\ (a,\alpha)\in I^p$, make up a basis of
$\Omega_p^n$. Set $I=\coprod_p I^p$.

An arbitrary element
$x_p\in\Omega_p^n\otimes X^p$ takes form
$$
x_p=\sum_{(a,\alpha)\in I^p}\omega_{a,\alpha}\otimes x_{a,\alpha}
$$
with $x_{a,\alpha}\in X^p$. This way we get a mapping
\begin{equation}
\label{mapping}
x\mapsto\{ x_{a,\alpha}\}_{(a,\alpha)\in I}
\end{equation}
from $\Omega^n(X)$ to the set of collections
\begin{equation}
\{ x_{a,\alpha}\}_{(a,\alpha)\in I},\ x_{a,\alpha}\in X^p\
\mbox{for}\ (a,\alpha)\in I^p
\label{collect}
\end{equation}

\subsubsection{} Let us introduce two operations on $I$.

For $(a,\alpha)\in I^p$ denote by
$\eta_i(a,\alpha)\in I^{p+1},\
i=1,\ldots,p+1$ the element
obtained from $(a,\alpha)$ by inserting $(0,0)$ on the place $i$.

Further, denote by $\zeta_i(a,\alpha)\in I^{p-1},\ i=1,\ldots,p-1$,
the pair $(a',\alpha')$ with $a'=(\ldots,a_i+a_{i+1},\ldots)$ and
$\alpha'=(\ldots,\alpha_i+\alpha_{i+1},\ldots)$--- this operation is defined
only if $\alpha_i+\alpha_{i+1}\leq 1$.

\subsubsection{} Let
$ x_p=\sum_{(a,\alpha)\in I^p}\omega_{a,\alpha}\otimes x_{a,\alpha}\in
\Omega^{n}_p\otimes X^p$
and
$ x_{p+1}=\sum_{(b,\beta)\in I^{p+1}}\omega_{b,\beta}\otimes x_{b,\beta}\in
\Omega^{n}_{p+1}\otimes X^{p+1}$.
The condition~(\ref{d-eq}) is equivalent to the following
formulas~(\ref{d_i-eq}) and~(\ref{d_0-eq}):

\begin{equation}
x_{\eta_i(a,\alpha)}=\delta^ix_{a,\alpha},\ \mbox{for }i\geq 1
\label{d_i-eq}
\end{equation}

\begin{equation}
 \delta^0(x_{a,\alpha})=
\sum\begin{Sb}
      e_0;\ e_i\leq a_i\\
      \epsilon_i\leq\alpha_i\\
      \beta_1:=\sum\epsilon_i\leq 1
    \end{Sb}
(-1)^{e+\beta_1+\sum_{i\geq j\geq 2}\epsilon_i\beta_j}{b_1!\over
{e_0!\cdots e_p!}} x_{b,\beta}
\label{d_0-eq}
\end{equation}
where the (big) sum is taken over non-negative $e_i,\epsilon_j$ satisfying
the conditions written and
$b_1=\sum_{i=0}^pe_i,
\ e=\sum_{i=1}^pe_i,\ b_{i+1}=a_i-e_i,\
\beta_1=\sum\epsilon_i,
\beta_{i+1}=\alpha_i-\epsilon_i.$

The condition~(\ref{s-eq}) is equivalent to the following
formulas~(\ref{s_i-eq}) and~(\ref{s_0-eq}):
\begin{equation}
\sigma^i(x_{b,\beta})=\begin{cases}
0 &\text{ if } \beta(i)=\beta(i+1)=1\\
\binom{b_i+b_{i+1}}{b_i}x_{\zeta_i(b,\beta)}& \text{ otherwise }
\end{cases}
\label{s_i-eq}
\end{equation}
for $i\geq 1$;
\begin{equation}
\sigma^0(x_{b,\beta})=\begin{cases}
x_{a,\alpha}& \text{ if } (b,\beta)=\eta_1(a,\alpha)\\
0 & \text{ otherwise.}
\end{cases}
\label{s_0-eq}
\end{equation}

Let us denote $\bT(X)$ the set of all collections ~(\ref{collect}) satisfying
{}~(\ref{d_i-eq})---~(\ref{s_0-eq}).

Let us call a collection ~(\ref{collect}) {\em locally finite} if
for every $p$ the set $\{ (a,\alpha)\in I^p|x_{a,\alpha}\neq 0\}$ is finite.
Let us denote by $\bT^{lf}(X)\subset\bT(X)$ the subset consisting of all
locally finite collections.

The above argument proves

\subsubsection{}
\begin{lem}{} The mapping ~(\ref{mapping}) defines an isomorphism
$$
\rho:\Omega^n(X)\iso\bT^{lf}(X)
$$
\end{lem} $\Box$

Elements $x_{a,\alpha},\ a,\alpha\in I$ are coordinates of $x\in\Omega^n(X)$.
Now our strategy will be: using relations ~(\ref{d_i-eq})---~(\ref{s_0-eq})
to express these coordinates in terms of smaller subsets of coordinates.

\subsection{} An element $(a,\alpha)\in I^p$ is called {\em reduced} if
none of $(a_i,\alpha_i),\ i=1,\ldots, p,$ is equal to $(0,0)$.

An element $(a,\alpha)\in I^p$ is called {\em special} if it
is reduced and $(a_1,\alpha_1)=(1,0)$.

An element $(a,\alpha)\in I^p$ is called {\em d-free} if it is  reduced
and not special.

The set of all d-free elements in $I^p$ will be denoted by ${\cal F}^p$.
We set ${\cal F}=\cup{\cal F}^p\subset I$.

Let us denote by $\bT_{\CF}(X)$ the set of all collections
\begin{equation}
\label{coll-f}
\{ x_{b,\beta}\}_{(b,\beta)\in\CF}, \mbox{ where } x_{b,\beta}\in X^p
\mbox{ for } (b,\beta)\in\CF^p
\end{equation}
satisfying following conditions:
\begin{align}
\sigma_0(x_{b,\beta})&=0; \nonumber\\
\sigma_i(x_{b,\beta})&=0 \text{ if } \beta(i)=\beta(i+1)=1\\
\sigma_i(x_{b,\beta})&=\binom{b_i+b_{i+1}}{b_i}x_{\zeta_i(b,\beta)}
\text{ otherwise } \nonumber
\label{props}
\end{align}

\subsubsection{}
\begin{lem}{} The natural projection defines an isomorphism
$$
\pi_1: \bT(X)\iso\bT_{\CF}(X)
$$
\end{lem}

{\bf Proof.} The formula~(\ref{d_i-eq}) allows one to express $x_{b,\beta}$ for
non-reduced $(b,\beta)\in I^{p+1}$ through
$x_{a,\alpha},(a,\alpha)\in I^p$. Furthermore, the formula~(\ref{d_0-eq})
allows one to calculate $x_{b,\beta}$ for special $(b,\beta)\in I^{p+1}$
by induction as follows. Endow $I^p$ with the partial order in which
$(a,\alpha)\geq (a',\alpha')$ iff $a_i\geq a_i'$ and $\alpha_i=\alpha'_i$
for all $i$. We determine the value of $x_{b,\beta}$ for special
$(b,\beta)\in I^{p+1}$
by induction on $(a,\alpha)$ such that $b=a\cup 1_1,\ \beta=\alpha
\cup 0_1$ in the obvious notation. For this one should consider the
equation~(\ref{d_0-eq}) and
see that all special summands in the right hand side except of $x_{b,\beta}$
correspond to smaller values of $(a,\alpha)$.

This immediately implies that
the map which takes $x=\{x_{a,\alpha}\}_{(a,\alpha)\in I^{\cdot}}\in\Th(X)^n$
to the collection
$\{x_{b,\beta}\}_{(b,\beta)\in {\cal F}^{\cdot}}$, is injective.

To prove bijectivity, we proceed by induction on $p$. Suppose that, apart
from $x_{b,\beta}$ with d-free $(b,\beta)$, all elements $x_{a,\alpha}$ with
$({a,\alpha})\in I^i,\ i\leq p$, are constructed and satisfy the
equations~(\ref{d_i-eq})--(\ref{s_0-eq}). In order to make the next
induction step, we have to check the following claims:

1) If $(b,\beta)=\eta_i(a,\alpha)=\eta_j(a',\alpha')$
then $\delta^ix_{a,\alpha}=\delta^jx_{a',\alpha'}.$

2) The condition~(\ref{d_0-eq}) is satisfied for any $(a,\alpha)\in I^p$
(and not only for reduced $(a,\alpha)$).

3) The conditions~(\ref{s_i-eq}) and~(\ref{s_0-eq}) are satisfied
for any $(b,\beta)\in I^{p+1}$ (and not only for d-free $(b,\beta)$).

These claims can be checked by a direct calculation using standard
identities for the morphisms in the category $\Delta$. $\Box$

\subsection{} An element $(a,\alpha)\in I$ is called {\em basic} if

--- $(a_1,\alpha_1)=$ either $(2,0)$ or $(0,1)$;

--- $(a_i,\alpha_i)=$ either $(1,0)$ or $(0,1)$ for $i>1$.

If $(a_1,\alpha_1)=(2,0)$ then $(a,\alpha)$ is called {\em  basic
element of the first kind}, otherwise --- {\em basic element of the
second kind}.

We denote $\CB$ the set of all basic elements; we define
$\CB^p=\CB\cap I^p$. Evidently $\CB^p\subset\CF^p\subset I^p$.

\subsubsection{}
\begin{lem}{} Suppose we are given $x=(x_{a,\alpha})_{(a,\alpha)\in I}\in
\bT(X)$. The elements $x_{a,\alpha},\ (a,\alpha)\in\CB$ satisfy the
following relations:
\begin{equation}
\label{rel-1}
\sigma^0(x_{a,\alpha})=0
\end{equation}

\begin{equation}
\sigma^i(x_{a,\alpha})=0 \text{ if } \alpha_i=\alpha_{i+1}=1
\end{equation}

\begin{equation}
\sigma^i(x_{a,\alpha})=\sigma^i(x_{a',\alpha'})
\end{equation}
for $i>1$ where
$$ a=(\ldots,0,1,\ldots), a'=(\ldots,1,0,\ldots)$$
$$ \alpha=(\ldots,1,0,\ldots), \alpha'=(\ldots,0,1,\ldots).$$
(switching the places $i,i+1$);

\begin{equation}
\sigma^1(x_{a,\alpha})=(\sigma^1)^2(x_{a',\alpha'})
\label{rel-3}
\end{equation}
where
$$ a=(2,0,\ldots),\ a'=(0,1,1,\ldots)$$
$$ \alpha=(0,1,\ldots),\ \alpha'=(1,0,0,\ldots).$$
\end{lem}

{\bf Proof.} Direct. $\Box$

\subsubsection{} We denote by $\bT_{\CB}(X)$
the set of all collections $\{ x_{a,\alpha}\}_{(a,\alpha)\in \CB}$
satisfying the relations ~(\ref{rel-1}) -- ~(\ref{rel-3}).

Thus we get a map
\begin{equation}
\label{map-pi2}
\pi_2:\bT_{\CF}(X)\lra \bT_{\CB}(X)
\end{equation}

\subsubsection{}
\begin{lem}{} The map ~(\ref{map-pi2}) is an isomorphism.
\end{lem}

{\bf Proof.} One can easily see that any $d$-free element may be obtained
from an element from $\CB$ by applying operations $\zeta_i$. This proves
injectivity of ~(\ref{map-pi2}). The proof of surjectivity is standard.
$\Box$

Summing up, we get a sequence of natural maps
\begin{equation}
\label{map-pi}
\Omega^n(X)\iso \bT^{lf}(X)\hra\bT(X)\iso\bT_{\CF}(X)\iso\bT_{\CB}(X)
\end{equation}

\subsection{} Let us call a collection
$\{x_{a,\alpha}\}_{(a,\alpha)\in{\cal B}^{\cdot}}$ {\em locally finite} if
$$
\forall p\in{\Bbb N}\ \exists m\in{\Bbb N}\ \forall(a,\alpha)\in
{\cal B}^q_{m'}
\text{ with }m'\geq m\ \forall f:[q]\ra[p]\in\Delta\text{  one has }
f(x_{a,\alpha})=0.
$$
Let us denote by $\bT^{lf}_{\CB}(X)\subset\bT_{\CB}(X)$ the subspace of all
locally finite collections.

\subsubsection{}
\begin{lem}{} The map ~(\ref{map-pi}) induces an isomorphism
$$
\pi:\Omega^n(X)\iso \bT^{lf}_{\CB}(X)
$$
\end{lem}

{\bf Proof.} Direct check. $\Box$

\subsection{} Set $\CB^p_m=\{ (a,\alpha)\in\CB^p|\sum a_i=m\}$;
$\CB_m=\cup_p\ \CB_m^p$. By definition,
$$
\CB_m=\CB_m^{m+n-1}\cup\CB_m^{m+n},
$$
the first (resp., second) subset consisting of all elements of the
first (resp., second) kind.

Let us introduce the following total order on $\CB_m$:
$(a,\alpha)>(a',\alpha')$ iff $a>a'$ in the lexicographic order.

Let us denote by $\bT_{\CB,m}(X)$ the set of all
collections $\{ x_{a,\alpha}\}_{(a,\alpha)\in \CB_m}$  satisfying the
relations ~(\ref{rel-1}) -- ~(\ref{rel-3}).

Given $(b,\beta)\in\CB_m$, denote by $\bT_{\leq (b,\beta)}(X)$ the
space
of all collections $\{ x_{a,\alpha}\}_{(a,\alpha)\in \CB_m,\ (a,\alpha)
\leq (b,\beta)}$.

We have obvious maps
$$
\bT_{\leq (b,\beta)}(X)\lra\bT_{\leq (a,\alpha)}(X)
$$
for $(a,\alpha)\leq (b,\beta)$.

\subsubsection{}
\label{inverse}
\begin{lem}{} We have
$$
\bT_{\CB,m}(X)=\invlim\ \bT_{\leq (b,\beta)}(X)
$$
the inverse limit over $\CB_m$.
\end{lem}

{\bf Proof.} Obvious. $\Box$

\subsection{}
\label{step}
\begin{lem}{}
Let $(a,\alpha)\in{\cal B}_m$. Let $(b,\beta)\in{\cal B}^l_m$ be the first
element such that $(b,\beta)>(a,\alpha)$.
Then the map
$$\bT_{\leq(b,\beta)}(X)\ra\bT_{\leq(a,\alpha)}(X)$$
is surjective and its kernel is isomorphic to a direct summand of $X^l$.
\end{lem}

{\bf Proof}
In order to lift an element of $\bT_{a,\alpha}(X)$, we have to find an
element $x_{b,\beta}$ having prescribed values for some of
 $\sigma^i(x_{b,\beta})$.

1. Existence of the lifting:

One has the following conditions on $\sigma^i(x_{b,\beta})$:

   (0) $\sigma^0(x_{b,\beta})=0$ --- always.

   (a) if $\beta_i=\beta_{i+1}=1$. Then one has $\sigma^i(x_{b,\beta})=0$.

   (b) if $(\beta_i,\beta_{i+1})=(0,1)$ and $i>1$. Then the condition is
$\sigma^i(x_{b,\beta})=\sigma^i(x_{b',\beta'})$ where the pair $(b',\beta')$
is obtained from $(b,\beta)$ by switching the places $(i,i+1)$.

   (c) if $(\beta_1,\beta_2)=(0,1)$. Then the condition is
$\sigma^1(x_{b,\beta})=(\sigma^1)^2(x_{a,\alpha})$ where
$(\alpha_1,\alpha_2,\alpha_3)=(1,0,0),\alpha_i=\beta_{i-1}\text{ for }
i>3$ and $a_i$
are defined uniquely by $\alpha_i$.

Let $I\subseteq[0,l-1]$ be the set of indices $i$ such that
$\sigma^i(x_{b,\beta})$ is defined by the conditions (0)---(c) above.
Let $y^i\in A^{l-1},\ i\in I,$ be the right hand sides of the
conditions (0)---(c) so that they take form
$$ \sigma^i(x_{b,\beta})=y^i,\ i\in I.$$

The first observation is that for any couple $i<j$ in $I$ one has
\begin{equation}
\sigma^{j-1}y^i=\sigma^iy^j.
\label{compatibility}
\end{equation}
 This can be  checked easily by an
explicit calculation.
Now consider two different cases:

(1st case) $I$ does not coincide with $[0,l-1]$.
In this case \Lem{kan} below asserts the existence of a solution
for the system of equations for $x_{b,\beta}$.

(2nd case) $I=[0,l-1]$. This is possible only in two cases:

--- $\beta_i=1$ for all $i$. Then $y^i=0$ for all $i$.

--- $\beta_1=0;\beta_i=1$ for $i>1$. Then $y^i=0$ for $i\not=1$ and
the condition~(\ref{compatibility}) gives that $y^1\in N^{l-1}X$.

In both cases \Lem{more} assures the existence of a solution.

2. The kernel of the map.

The kernel of the map in question has always form
$$ N^p_I(X):=\{a\in X^p| \sigma^i(a)=0\text{ for all }i\in I\}$$
for some subset $I\subseteq[0,p]$. According to~\Lem{kan2} below, this is
a direct summand of $X^p$. $\Box$

\subsection{} Let $\cal M$ be the set of all non-decreasing functions
$f:{\Bbb N}\ra{\Bbb N}$ equipped with a partial order:
$$
f\geq g\text{ iff }f(n)\geq g(n)\text{ for each }n.
$$

\subsubsection{}
\begin{defn}{growth}
Given $f\in{\cal M}$, an element $x\in\bT(X)$ {\em has growth $\leq f$}
if
$\pi_1(x)=\{x_{b,\beta}\}_{(b,\beta)\in{\cal F}}$
satisfies the property
$$
x_{b,\beta}=0 \text{ if } (b,\beta)\in{\cal F}^p\text{ and }
\sum b_i>f(p).
$$
Denote by $\bT^f(X)\subset \bT(X)$ the subspace of all elements of growth
$\leq f$.
\end{defn}

We have $\bT^f(X)\subseteq \bT(X)$. Moreover,

\subsubsection{}
\label{directlim}
\begin{lem}{} We have
$$
\bT^{lf}(X)=\dirlim\ \bT^f(X)
$$
the limit taken over $f\in\CM$.
\end{lem}

{\bf Proof.} Obvious.

\subsection{}
For $X\in\Delta\Mod(A)$ and $d\in{\Bbb N}\cup\{ -1,\infty\}$ define
$X_{>d}\in\Delta\Mod(A)$ as follows. Set $X_{>-1}=X;\ X_{>\infty}=0$.
For $d\in \Bbb N$ set
\begin{equation}
X_{>d}^i=\{x\in X^i|f(x)=0\ \forall f:[i]\ra[d]\in\Delta\}
\end{equation}

We have obviously $X_{>d}^i=0\text{ for }i\leq d$, and $X_{>d}^{d+1}=N^{d+1}X$
where $N^iX$ denotes the normalization (see ~\ref{normal} below).

\subsubsection{} Given $f\in\CM$, define a function
$f^{\circ}: \Bbb N\lra \Bbb N\cup\{ -1,\infty\}$ by a formula
$$
f^{\circ}(m)=p\text{ iff } f(p)<m\leq f(p+1).
$$

For $f\in\CM$ let us denote
$$
\bT_{\CB}^f(X)=\pi_2\pi_1(\bT^f(X))
$$

\subsubsection{}
\begin{lem}{} We have
$$
\bT_{\CB}(X)=\prod_{m=0}^{\infty}\bT_{\CB,m}(X)
$$
\end{lem}
{\bf Proof.} Evident: the equations~(\ref{rel-1}--\ref{rel-3})
are homogeneous on $m=\sum a_i$.
 $\Box$

\subsubsection{}
\begin{lem}{product} For any $f\in{\cal M}$ one has
\begin{equation}
\label{eq-prod}
\bT^f_{\CB}(X)=\prod_{m=0}^{\infty}\bT_{\CB,m}(X_{>f^{\circ}(m)})
\end{equation}
\end{lem}

\begin{pf} An element $x\in\bT(X)$ has growth $\leq f$ iff
$\pi_1(x)=\{x_{b,\beta}\}_{(b,\beta)\in{\cal F}}$ satisfies the property
$$
x_{b,\beta}=0 \text{ if } (b,\beta)\in{\cal F}^p\text{ and }
p\leq f^{\circ}(\sum b_i).
$$
Then the formulas~(\ref{props}) give immediately the result.
\end{pf}

Recall that $X$ is called {\em finite dimensional} if there exists
$d\in\Bbb N$
such that $N^i(X)=0$ for $i>d$ (we say that $\dim(X)\leq d$.

If this is the case then $X_{>d}=0$, hence the product in ~(\ref{eq-prod})
is finite.

Combining all the previous results together, we get

\subsection{}
\begin{thm}{} For every $n\in\Bbb N$ the $A$-module $\Omega^n(X)$ may be
obtained from $X$ applying the following operations (naturally in $X$):

(1) taking modules $\bT_{\CB,m}(X)$ which have a natural finite filtration
with graded factors isomorphic to direct summands of modules
$X^l,\ l\in\Bbb N$;

(2) taking direct products over $\Bbb N$; if $\dim(X)<\infty$ then
products are finite;

(3) passing to a filtered direct limit.

\end{thm}

{\bf Proof.} Follows immediately from ~\ref{directlim}, ~\ref{product},
{}~\ref{step} and ~\ref{inverse}. $\Box$

We have the following easy corollaries which are the main properties
of the functor $\Omega$.

\subsection{}
\label{exact-6}
\begin{cor}{} The functor $X\mapsto\Omega(X)$ is exact. $\Box$
\end{cor}

\subsection{}
\label{base-6}
\begin{cor}{} Suppose $X\in\Delta\Mod(A)$ is finite dimensional. Let $A'$ be
a (commutative) $A$-algebra. Then the natural base change map
$$
\Omega(X)\otimes_AA'\lra\Omega(X\otimes_AA')
$$
is an isomorphism. $\Box$
\end{cor}

\subsection{}
\label{flat-6}
\begin{cor}{}
Let $X\in\Delta\Mod(A)$ and suppose that either $A$ is noetherian
or $X$ is finite dimensional.

If all
$X^i$ are flat $A$-modules so are all $\Omega^n(X)$. $\Box$
\end{cor}

\subsection{}
\label{qc}
\begin{cor}{} Let $S$ be topological space endowed with a sheaf $\CO_S$
of commutative $\Bbb Q$-algebras. Let $X\in\Delta\Mod({\CO}_S)$ be
finite. If all
$X^i$ are quasi-coherent ${\CO}_S$-modules so are all $\Omega^n(X)$.
$\Box$
\end{cor}

In the remaining part of this Section we will prove some facts about
cosimplicial modules needed above. Most of them are more or less standard.

In fact, mostly we will discuss an explicit form of Dold-Puppe correspondence.

\subsection{}
\label{normal}
Let $X$ be a cosimplicial $A$-module. For any $i\geq 0$ define
$$ N^i(X)=\{x\in X^i|\sigma^j(x)=0\text{ for all }j\}.$$

Define a subcategory $\Lambda$ in $\Delta$ as follows. $\Lambda$ has
the same objects as $\Delta$; The set of morphisms of $\Lambda$ is generated
by the faces $\delta^i:[n]\ra[n+1]$ with $i=0,\ldots, n$
(that is: (a) only faces appear; (b) the last face $\delta^{n+1}:[n]\ra[n+1]$
disappear).

Define {\em a shift functor} $  s:\Delta\ra\Delta$
by the formulas
$$s[n]=[n+1],\ s(\delta^i)=\delta^{i+1},\ s(\sigma^i)=\sigma^{i+1}.$$

\subsubsection{}
\begin{prop}{decomposition}
For any cosimplicial $A$-module $X$ one has
$$
X^n=\bigoplus_{m\geq 0}\bigoplus_{f:[m]\ra[n]\in\Lambda}f(N^m(X)).
$$
\end{prop}
\begin{pf}

We will prove the claim by induction.

For $n=0$ the claim is trivial. Suppose it is true for degrees $< n$
and for all cosimplicial modules $X$. Apply this to the shift $Y=Xs$
of $X$. We have by the inductive hypothesis
$$ X^n=Y^{n-1}=\bigoplus_{m\geq 0}
\bigoplus_{f:[m]\ra[n-1]\in{\Lambda}}f(M^m)$$
where $M^i=N^i(Y)$ is the normalization of $Y=Xs$.
An element $x\in X^{m+1}$ belongs to $M^m$ iff $\sigma^i(x)=0$ for
$i>0$. Write $x=y+\delta^0\sigma^0(x)$ where $y=x-\delta^0\sigma^0(x)$.
One checks that $y\in N^{m+1}(X)$ and (of course) $\sigma^0(x)\in X^m$.
Thus by induction (note that $m\leq n-1$)
$$\sigma^0(x)=\sum_k\sum_{g:[k]\ra[m]\in\Lambda}g(z_g)$$
with $z_g\in N^k(X)$.
If $x'\in M^m$ is the element corresponding to $x\in X^{m+1}$ then
$f(x')=(sf)(x)$ and therefore
$$ f(x')=(sf)(x)=(sf)(y)+(sf)\delta^0\sum g(z_g)=(sf)(y)+\sum\delta^0
fg(z_g)$$
which proves that
$$ X^n=\sum_m\sum_{f:[m]\ra[n]\in\Lambda}f(N^m(X)).$$
Let us finally prove the uniqueness of the presentation
of an element $x\in X^n$ into sum

\begin{equation}
 x=\sum_m\sum_{f:[m]\ra[n]\in\Lambda}f(x_f)
\label{presentation}
\end{equation}
with $x_f\in N^m(X)$.

For any $f:[m]\ra[n]\in\Lambda$ define a left-inverse $f^l:[n]\ra[m]$
as follows: if $f=\delta^{i_1}\cdots\delta^{i_{n-m}}$ with
$i_1>\ldots>i_{n-m}$ then
$$ f^l=\sigma^{i_{n-m}}\cdots\sigma^{i_1}.$$

We will check the uniqueness of the elements $x_f$ in the
presentation~(\ref{presentation}) by induction on $m$.

Suppose that $x_f$ are defined uniquely for all $f:[m']\ra[n]$ with $m'<m$
(this is true, say, for $m=0$).
Define an order on the set $\Hom_{\Lambda}([m],[n])$ saying that
$\delta^{i_1}\cdots\delta^{i_{n-m}}\geq\delta^{j_1}\cdots\delta^{j_{n-m}}$
iff $(i_1,\ldots,i_{n-m})\geq(j_1,\ldots,j_{n-m})$ in the lexicographic order.
One immediately checks the following

\begin{lem}{}
If $f>g$ then $f^lg\not=\id_{[m]}$.
\end{lem}

By the inductive hypothesis we can suppose that $x_f=0$ for $f:[m']\ra[n]$
with $m'<m$. Then
\begin{equation}
f^l(g(x_g))=\begin{cases}
x_g&\text{ if }f=g\\
0& \text{ if } f>g
\end{cases}
\end{equation}
Thus, the transition matrix expressing the values of $f^l(x)$ for different $f$
through $x_g$ is upper-triangular and hence invertible.

The proposition is proven.
\end{pf}

It is very convenient to rewrite the statement of~\Prop{decomposition} as
follows.

Let ${\Bbb Q}\Lambda_{mn}$ be the rational vector space spanned by
the set $\Hom_{\Lambda}([m],[n])$. Then one has

\begin{cor}{}
One has a canonical on $X\in\Delta\Mod(A)$ isomorphism
$$ X^n=\bigoplus_{m\geq 0}{\Bbb Q}\Lambda_{mn}\otimes_{\Bbb Q}N^m(X).$$
\end{cor}

\subsection{}
\begin{lem}{kan}
Let $I\subset[0,n]$ be a proper subset and let $y_i\in X^{n-1},\ i\in I$
be given.
Then the system of equations
$$ \sigma^i(x)=y^i,\ i\in I,$$
has a solution if (and only if) $y^i$ satisfy the compatibilities
$$ \sigma^{j-1}y^i=\sigma^iy^j\text{ for }i<j\in I.$$
\end{lem}
\begin{pf}
Consider the functor
$$ i:\Delta^0\ra\Delta $$
defined by the formulas
$$ i([n])=[n+1],\ i(\delta^i)=\sigma^i,\ i(\sigma^i)=\delta^{i+1}$$
--- see Gabriel-Zisman's functor II, ~\cite{gz}, 3.1.1).

If $X$ is a cosimplicial module then $Xi$ is a simplicial module.
The property of $X$ we have to prove just means that $Xi$ is a
Kan simplicial set. This is well-known to be always true
(see, e.g., ~\cite{q2}, Prop. II.3.1).
\end{pf}

\subsection{}
\begin{lem}{more}For any system $y^i\in N^{n-1}X,\ i\in[0,n-1]$
there exists an element $x\in X^n$ such that $\sigma^ix=y^i$ for
all $i\in[0,n-1]$.
\end{lem}
\begin{pf}
We will look for $x=\sum\delta^i(z^i)$ with $z^i\in N^{n-1}X$.
Then the conditions on $x$
are expressed by the equations $y^i=\sigma^i(x)=z^i+z^{i+1}$
which are clearly solvable.
\end{pf}

\subsection{}
{\em Notation.} For $I\subseteq[0,n]$ denote
$$ N^i_I(X)=\{x\in X^i|\sigma^j(x)=0\text{ for all }j\in I\}.$$

\subsubsection{}
\begin{lem}{kan2}
Let $I\subseteq[0,n]$. There exist a collection of vector subspaces
${\Bbb Q}\Lambda^I_{m,n}\subseteq{\Bbb Q}\Lambda_{m,n}$
such that
$$
N^n_I(X)=\bigoplus_m{\Bbb Q}\Lambda^I_{m,n}\otimes_{\Bbb Q}N^m(X)
$$
\end{lem}

Here is a more general statement.

\subsubsection{}
\begin{lem}{kan3}
Fix $d\leq n\in{\Bbb N}$ and a subset $\Phi\subseteq\Hom_{\Delta}([n],[d])$.
Define
$$X^n_{\Phi}=\{x\in X^n|f(a)=0\text{ for all }f\in\Phi\}.$$

 There exists a collection of vector subspaces
${\Bbb Q}\Lambda^{\Phi}_{m,n}\subseteq{\Bbb Q}\Lambda_{m,n}$
such that
$$ X^n_{\Phi}=\bigoplus_m{\Bbb Q}\Lambda^{\Phi}_{m,n}\otimes_{\Bbb Q}N^m(X).$$
\end{lem}

\begin{pf} of~\ref{kan3}.
Let $t:[n]\ra[k]\in\Delta$.
The condition $tx=0$ for $x=\sum_ff(x_f)$ takes form
$$0=tx=\sum_ftf(x_f)=\sum_gg(\sum_{f: tf=g}x_f)$$
which is equivalent to the system of equations
$$ \sum_{f: tf=g}x_f=0$$
numbered by $g\in\Mor\Delta$.

The lemma immediately follows from this observation.
\end{pf}

\section{Higher Kodaira-Spencer morphisms}
\label{kodaira}

\subsection{}
\label{laawm}
(Cf. \cite{bs}, 1.2.) Let $S$ be a differentially smooth formal
scheme (for example, a smooth scheme).   Let $\pi: X\lra S$
be a smooth separated map of formal schemes, ~\ref{smooth}, for example
a smooth morphism of usual schemes.
We have the exact sequence ~(\ref{seq}):
$$
0\lra \CT_{X/S}\lra \CT_X\overset{\epsilon_{\CT}}{\lra} \pi^*\CT_S\lra 0
$$
The first embedding makes $\CT_{X/S}$ a Lie algebroid over $X$. Note that
the sheaf $\pi^*\CT_S$ is not a sheaf of Lie algebras.

Let $\pi^{-1}$ denotes the functor of set-theoretical inverse image, so that
$\pi^*=\CO_X\otimes_{\pi^{-1}\CO_S}\pi^{-1}$. The subsheaf
$\pi^{-1}\CT_S\subset \pi^*\CT_S$ is a $\pi^{-1}\CO_S$-Lie algebra.
Set $\CT_{\pi}:=\epsilon_{\CT}^{-1}(\pi^{-1}\CT_S)$ - this is a sheaf
of $k$-Lie algebras and $\pi^{-1}\CO_S$-modules (consisting of vector fields
with the constant projection to $S$ along fibers of $\pi$).
We have an exact sequence of $k$-Lie algebras and $\pi^{-1}\CO_S$-modules
$$
0\lra\CT_{X/S}\lra \CT_{\pi}\lra \pi^{-1}\CT_S\lra 0
$$

Let $\epsilon: \CA\lra \CT_X$ be a transitive Lie algebroid over $X$.
We set $\CA_{/S}:=\epsilon^{-1}(\CT_{X/S})\subset\CA_{\pi}:=
\epsilon^{-1}(\CT_{\pi})\lra \CA$. These are subsheaves of Lie
algebras. $\CA_{/S}$ is a subsheaf of $\CO_X$-modules, and a Lie algebroid
over $X$ included into an exact sequence
$$
0\lra \CA_{(0)}\lra \CA_{/S}\lra \CT_{X/S}\lra 0
$$
$\CA_{\pi}$ is a subsheaf of $\pi^{-1}\CO_S$-modules. We have an exact
sequence of $k$-Lie algebras and $\pi^{-1}\CO_S$-modules
\begin{equation}
\label{a}
0\lra\CA_{/S}\lra\CA_{\pi}\lra \pi^{-1}\CT_S\lra 0
\end{equation}

\subsection{}
\label{dglie-fund} Pick a finite affine open covering $\CU$ of $X$.
We will suppose that $\CA$ is a locally free $\CO_X$-module of finite
type, hence, so is $\CA_{/S}$.
Let us apply the functor $\pi_*\CHC(\CU,\cdot)$ to the exact sequence
{}~(\ref{a}). We have $R^1\pi_*\CHC^i(\CU,\CA_{/S})=0$; so  we get
an exact sequence of cosimplicial $\CO_S$-modules
$$
0\lra \pi_*\CHC(\CU,\CA_{/S})\lra \pi_*\CHC(\CU,\CA_{\pi})
\lra \pi_*\CHC(\CU,\pi^{-1}\CT_S)\lra 0
$$
Next, applying the Thom-Sullivan functor $\Omega$ we get an exact sequence
of complexes
\begin{equation}
\label{ex-derham}
0\lra \pi^{Lie}_{*,\CU}(\CA_{/S})\lra \Omega(\pi_*\CHC(\CU,\CA_{\pi}))
\lra \Omega(\pi_*\CHC(\CU,\pi^{-1}\CT_S))\lra 0
\end{equation}
Note that $\Omega(\pi_*\CHC(\CU,\CA))$ is naturally a dg $k$-Lie algebra.
We have a canonical adjunction map
$$
\CT_S\lra \Omega(\pi_*\CHC(\CU,\pi^{-1}\CT_S))
$$
so, taking the pullback of ~(\ref{ex-derham}) we get an exact sequence
\begin{equation}
\label{fund}
0\lra \pi^{Lie}_{*,\CU}(\CA_{/S})\lra \CA^{\pi}_{\CU}\lra \CT_S\lra 0
\end{equation}
By definition,
$$
\CA^{\pi}_{\CU}=\Omega(\pi_*\CHC(\CU,\CA_{\pi}))
\times_{\Omega(\pi_*\CHC(\CU,\pi^{-1}\CT_S))}\CT_S,
$$
and this sheaf inherits the structure of a dg $k$-Lie algebra and
$\CO_S$-module from $\Omega(\pi_*\CHC(\CU,\CA_{\pi}))$ and $\CT_S$.
One sees that this way we get a structure of a transitive dg Lie algebroid on
$\CA^{\pi}_{\CU}$.

\subsection{}
\label{higher} Now we can apply to ~(\ref{fund}) the construction
{}~\ref{abstr-ksmaps}. We get the maps:
\begin{equation}
\label{ks-class}
\kappa^1:\CT_S\lra R^1\pi_*(\CA_{/S})
\end{equation}
--- the "classical" KS map;
\begin{equation}
\label{ks}
\kappa:\Diff_S\lra \CH_0^{Lie}(\pi_*^{Lie}(\CA_{/S}))
\end{equation}
and
\begin{equation}
\label{ks-leq-n}
\kappa^{\leq n}:\Diff^{\leq n}_S\lra \CF_n\CH^{Lie}_0(\pi_*^{Lie}(\CA_{/S}))
\end{equation}
satisfying the compatibilities ~\ref{main-thm}. These maps are called
{\bf higher Kodaira-Spencer morphisms.}

\subsection{Split case}
\label{dglie-fund-const} Suppose that $X=Y\times S$, and $\pi$ is a projection
to the second factor. In this case we have canonical embeddings
$$
\pi^{-1}\CT_S\hra\pi^*\CT_S\hra \CT_X
$$
By taking the pull-back of the exact sequence
$$
0\lra\CA_{(0)}\lra \CA\lra\CT_X\lra 0,
$$
we get an exact sequence
\begin{equation}
\label{a-const}
0\lra \CA_{(0)}\lra \bar{\CA}\lra\pi^{-1}\CT_S\lra 0
\end{equation}
Now we can repeat the construction of ~\ref{dglie-fund}, replacing
the sequence ~(\ref{a}) by ~(\ref{a-const}). This will
provide a dg Lie algebroid
\begin{equation}
\label{fund-const}
0\lra\CA_{(0)}\lra\bar{\CA}^{\pi}_{\CU}\lra\CT_S\lra 0
\end{equation}

Again we can apply the construction ~\ref{abstr-ksmaps},
and get the KS maps analogous to ~(\ref{ks-class}) -~(\ref{ks-leq-n}):
\begin{equation}
\label{ks-class0}
\kappa^1_{(0)}:\CT_S\lra R^1\pi_*(\CA_{(0)})
\end{equation}
--- the "classical" KS map;
\begin{equation}
\label{ks0}
\kappa_{(0)}:\Diff_S\lra \CH_0^{Lie}(\pi_*^{Lie}(\CA_{(0)}))
\end{equation}
and
\begin{equation}
\label{ks-leq-n0}
\kappa^{\leq n}_{(0)}:\Diff^{\leq n}_S\lra \CF_n\CH^{Lie}_0(\pi_*^{Lie}
(\CA_{(0)}))
\end{equation}
satisfying the same compatibilities.

\subsection{}
\label{isom}
\begin{thm}{} Suppose  that

--- either we are in the situation ~\ref{higher}, $\pi_*\CA_{/S}=0$,
and $\kappa^1$ is an isomorphism,

--- or we are in the situation ~\ref{dglie-fund-const}, $\pi_*\CA_{(0)}=0$
and $\kappa^1_{(0)}$ is an isomorphism.

Then all $\kappa^{\leq n}$ (resp., $\kappa^{\leq n}_{(0)}$) are isomorphisms.
\end{thm}

{\bf Proof.} Suppose for definiteness that we are in the first situation.
The claim is proved by induction on $n$, using commutative diagrams

\begin{center}
  \begin{picture}(20,4)
     \put(0,0){\makebox(1,1){$0$}}
     \put(2.5,0){\makebox(4,1){$\CF_{n-1}\CH^{Lie}_0(\pi^{Lie}_*(\CA_{/S}))$}}
     \put(8,0){\makebox(4,1){$\CF_{n}\CH^{Lie}_0(\pi^{Lie}_*(\CA_{/S}))$}}
     \put(13.5,0){\makebox(4,1){$S^n(R^1\pi_*(\CA_{/S}))$}}
     \put(19,0){\makebox(1,1){$0$}}

     \put(0,3){\makebox(1,1){$0$}}
     \put(2.5,3){\makebox(4,1){$\Diff^{\leq n-1}_S$}}
     \put(8,3){\makebox(4,1){$\Diff_S^{\leq n}$}}
     \put(13.5,3){\makebox(4,1){$S^n(\CT_S)$}}
     \put(19,3){\makebox(1,1){$0$}}

     \put(4.5,3){\vector(0,-1){2}}
     \put(10,3){\vector(0,-1){2}}
     \put(15.5,3){\vector(0,-1){2}}

    \put(1,0.5){\vector(1,0){1}}
    \put(7,0.5){\vector(1,0){1}}
    \put(12.5,0.5){\vector(1,0){1}}
    \put(17.8,0.5){\vector(1,0){1}}

    \put(1,3.5){\vector(1,0){2}}
    \put(5.5,3.5){\vector(1,0){3}}
    \put(11,3.5){\vector(1,0){3}}
    \put(17,3.5){\vector(1,0){2}}

    \put(2.5,1){\makebox(2,2){$\kappa^{\leq n-1}$}}
    \put(8,1){\makebox(2,2){$\kappa^{\leq n}$}}
    \put(16,1){\makebox(2,2){$(-1)^nS^n(\kappa^1)$}}

  \end{picture}
\end{center}

Note that our assumptions imply that
$$H^{-1}(F_nC(\pi^{Lie}_*(\CA_{/S}))/F_{n-1}C(\pi^{Lie}_*(\CA_{/S})))=0$$
and
$$H^{0}(F_nC(\pi^{Lie}_*(\CA_{/S}))/F_{n-1}C(\pi^{Lie}_*(\CA_{/S})))=
S^n(R^1\pi_*(\CA_{/S})).$$
$\Box$

\subsection{Deformations of schemes.}
\label{def-schemes} Set
$\CA=\CT_X$. Applying the previous construction, we get the maps
\begin{equation}
\label{ks-class-1}
\kappa^1: \CT_S\lra R^1\pi_*(\CT_{X/S})
\end{equation}
--- the classical Kodaira-Spencer map;
\begin{equation}
\label{ks-1}
\kappa: \Diff_S\lra \CH_0^{Lie}(\pi_*^{Lie}(\CT_{X/S}))
\end{equation}
and
\begin{equation}
\label{ks-leq-n-1}
\kappa^{\leq n}: \Diff_S^{\leq n}\lra \CF_n\CH_0^{Lie}(\pi_*^{Lie}(\CT_{X/S}))
\end{equation}
satisfying the compatibilities ~\ref{main-thm}.

\subsection{Deformations of $G$-torsors}
\label{def-tors} (a) Let $G$ be an algebraic group over $k$, $\fg=\Lie(G)$,
and $p:P\lra X$ a $G$-torsor over $X$. We define sheaves of Lie algebras
$\CA_P$ and $\fg_P$ as in ~\ref{pose}. $\CA_P$ is naturally a transitive Lie
algebroid over $X$, with $\CA_{P(0)}=\fg_P$.

Applying the previous construction, we get the maps analogous to
{}~(\ref{ks-class-1}) - ~(\ref{ks-leq-n-1}), with $\CT_{X/S}$
replaced by $\CA_{P/S}$, subject to the same compatibilities.

(b) Suppose in addition that $X=Y\times S$ as in ~\ref{dglie-fund-const}.
Then we can apply the construction of {\em loc. cit.} to $\CA=\CA_P$, and get
higher KS maps taking value in  $\CF_n\CH_0^{Lie}(\fg_P)$.

\section{Universal deformations}
\label{univers}

\subsection{} Let us fix a smooth scheme $X$, an algebraic group $G$ and
a $G$-torsor $P$ over $X$. Consider the following deformation
problems. To each problem we assign a sheaf of $k$-Lie algebras
$\CA_i,\ i=1,2,3$ over $X$.

{\bf Problem 1.} Flat deformations of $X$; $\CA_1=\CT_X$.

{\bf Problem 2.} Flat deformations of the pair $(X,P)$; $\CA_2=\CA_P$,
cf. \ref{def-tors}.

{\bf Problem 3.} Deformations of $P$ ($X$ being fixed); $\CA_3=\fg_P$.

Accordig to Grothendieck, to each problem corresponds a (2-)functor
of infinitesimal deformations
$$
F_i:\Artin_k\lra\Groupoids
$$
from the category of local artinian $k$-algebras with the residue
field $k$ to the (2-)category of groupoids.

In each case, $\CA_i$ is "a sheaf of infinitesimal automorphisms"
corresponding to $F_i$ (in the sense of ~\cite{sga1}, Exp.III, 5,
especially Cor. 5.2 for Problem 1; for the other problems the meaning is
analogous).

In particular, we have the Kodaira-Spencer-Grothendieck isomorphisms
$$
\ft_{\bF_i}\cong H^1(X,\CA_i).
$$
where
$$
\bF_i:\Artin_k\lra \Ens
$$
is the composition of $F_i$ and the connected components functor
$\pi_0:\Groupoids\lra \Ens$. Here for a functor
$$
F:\Artin_k\lra\Ens
$$
$\ft_F$ denotes the tangent space to $F$:
$$
\ft_F=F(\Spec(k[\epsilon]/(\epsilon^2)),
$$
cf. ~\cite{gr}, Exp. 195.

\subsection{} One can verify that in each case the conditions of
Schlessinger's Theorem 2.11, ~\cite{sch}, are fullfilled, and there exists
a versal formal deformation $\fS_i$.

\subsubsection{}
\begin{lem}{} Suppose that $H^0(X,\CA_i)=0$. Then $\fS_i$ is a universal
deformation, i.e. $\bF_i$ is prorepresentable.
\end{lem}

\begin{pf} This fact is presumably classic. We give a sketch of a proof for
completeness. We have $\fS_i=\Spf(R)$, and we have a canonical morphism
$x: h_R\lra F$ (we use notations of ~\cite{sch}). Consider the functor
$$
G:\Artin_R\lra\Ens
$$
from the category of local artinian $R$-algebras with the residue field $k$,
defined as $G(\alpha:R\lra A))=\Aut_{F_i(A)}(\alpha_*(x))$.

{\bf Claim 1.} $G$ is prorepresentable.

Indeed, one can check that the hypotheses of {\em loc.cit.}, 2.11 hold true
for $G$.

Let $\Spf(T)$ prorepresents $G$, where $T$ is a complete local $R$-algebra.

{\bf Claim 2.} The structure morphism $\phi: R\lra T$ is injective.

In fact, this morphism has a section since groups of automorphisms have
identity. Now, we have
$$
\bm_T/\bm^2_T+\bm_R=H^0(X,\CA_i)=0
$$
hence $R=T$, whence $\bF_i$ is prorepresentable by {\em loc.cit.} 3.12.
\end{pf}

\subsection{}
\label{complet}
\begin{thm}{} Suppose that $H^0(X,\CA_i)=0$; let $\fS=\Spf(R)$ be the base
of the universal formal deformation for Problem $i$.
Suppose that $\fS$ is formally smooth. Then we have a canonical
isomorphism
$$
R^*=H^{Lie}_0(R\Gamma^{Lie}(X,\CA_i))
$$
where $R^*$ denotes the space of continuous $k$-linear maps $R\lra k$
($k$ considered in the discrete topology).
\end{thm}

\begin{pf} Let us treat Problem 2 for definiteness (for other problems
the proof is the same).  Since $\fS$ is formally
smooth, $R$ is isomorphic to a power series algebra $k[[T_1,\ldots,T_n]]$,
hence $\fS$ is differentially smooth. Let
$\pi: \fX\lra\fS$ be the universal deformation and $\fP$ the
universal $G$-torseur over $\fX$. Applying ~\ref{def-tors}, we get the map
$$
\kappa:\Diff_{\fS}\lra \CH^{Lie}_0(\pi_*^{Lie}(\CA_{\fP/\fS}))
$$
Since $\fS$ is universal, the usual KS map
$\kappa^1:\CT_{\fS}\lra R^1\pi_*(\CA_{\fX/\fS})$ is isomorpism.
Note that  $\CA_{\fP/\fS}|_X\cong \CA_i$.
Since $H^0(X,\CA_i)$, we have $\pi_*(\CA_{\fX/\fS})=0$; hence by
{}~\ref{isom} $\kappa$ is isomorphism.

Now let us take the geometric fiber of $\kappa$ at the closed point $s:\Spec(k)
\hra \fS$. We have $\Diff(\fS)_{k(s)}\cong R^*$ by ~\ref{point} and
$$
\CH^{Lie}_0(\pi_*^{Lie}(\CA_{\fP/\fS}))_{k(s)}\cong
H^{Lie}_0(R\Gamma^{Lie}(X,\CA_i))
$$
by ~\ref{base}. The theorem follows.
\end{pf}



\end{document}